\documentclass[leqno]{article}
\usepackage[latin1]{inputenc}
\usepackage[francais,english]{babel}
\usepackage{amsmath}
\usepackage{amssymb}
\usepackage{epsfig}
\usepackage{graphicx}
\usepackage{color}
\usepackage{psfrag}

\newcommand{\cqfd}
{%
\mbox{}%
\nolinebreak%
\hfill%
\rule{2mm}{2mm}%
\medbreak%
\par%
}

\newtheorem{theorem}{Theorem}
\newtheorem{lemma}[theorem]{Lemma}
\newtheorem{proposition}[theorem]{Proposition}

\newtheorem{remark}[theorem]{Remark}
\newtheorem{claim}[theorem]{Claim}

\title{Distributed source identification for wave equations : an observer-based approach\thanks{This work was  partially supported  by the "Agence Nationale de la Recherche" (ANR),
     "Projet Jeunes Chercheurs" EPOQ2 number ANR-09-JCJC-0070.}}

\author{Marianne Chapouly and Mazyar Mirrahimi \thanks{INRIA Rocquencourt, Domaine de Voluceau, B.P. 105, 78153 Le Chesnay Cedex,  (marianne.chapouly@inria.fr, mazyar.mirrahimi@inria.fr)} }

\date{}

\begin{document}

\maketitle
\begin{abstract}
In this paper, we consider the 1D wave equation where the spatial domain is a bounded interval. Assuming the initial conditions to be known, we are here interested in identifying an unknown source term, while we take the Neumann derivative of the solution on one of the boundaries as the measurement output. Applying a back-and-forth iterative scheme and constructing well-chosen  observers, we retrieve the source term from the measurement output in the minimal observation time. 
\end{abstract}
\vspace{0.2 cm}


\vspace{0.2 cm} {\bf AMS subject classifications.} 
\section{Introduction}

\indent Blum and Auroux~\cite{2005-Auroux-Blum} proposed a new inversion algorithm for identifying the initial state of an observable system, based on the application of back-and-forth observers. Noting that we have only access to a measurement output on a fixed time interval $(0,T)$, the idea consists in proposing a first asymptotic observer for the system that will be applied in this time interval and a second one that will be applied to the system where the direction of time is reversed. These observers are then used iteratively to get a better estimate of the initial state after each back-and-forth iteration. If the two observer gains  are well-chosen, so that the whole back-and-forth procedure induces a contracting error dynamic, one can ensure the convergence of the estimator towards the initial state. 

More recently, Ramdani-et-al~\cite{2009-Ito-Ramdani-Tucsnak,2009-Ramdani-Tucsnak-Weiss} have considered, for the case of wave equations, the theoretical study of this problem applying techniques borrowed from semigroups theory.  Here, we consider a similar problem to~\cite{2009-Ramdani-Tucsnak-Weiss} for a wave equation, where the initial state is known (is zero) and we are rather interested in identifying an unknown source term $q(x)$. Let also $T>0$,  $\omega\in\mathbb{R}$ and let $q\in H^2(0,1)\cap H^1_0(0,1)$. We consider the following system
\begin{equation}
\label{eq1}
\left\{
\begin{array}{l}
u_{tt}-u_{xx}=q(x)\cos(\omega t),\,(t,x)\in(0,T)\times(0,1),\\
u(t,0)=u(t,1)=0,\,t\in(0,T),\\
u(0,x)=u_t(0,x)=0,\,x\in(0,1),\\
y(t)=u_x(t,0),\,t\in(0,T),
\end{array}
\right. 
\end{equation}
where $(u,u_t)$ represents the state of the system and $y$ is the output. The term $cos (\omega t)q(x)$, where $\omega$ is a fixed (known) frequency and $q(x)$ is unknown, is considered to be an external force which varies harmonically.\\
\indent For any $q\in H^2(0,1)\cap H^1_0(0,1)$ there exists a unique solution $u\in C^1([0,T];H^1_0(0,1))\cap C^2([0,T];L^2(0,1))$ to \eqref{eq1} and $u_x(t,0)\in H^1(0,T)$ (see \cite[Remark 2]{1995-Yamamoto}). It is moreover well-known from the work \cite{1995-Yamamoto} of Yamamoto that this problem is well-posed in the sense that one can retrieve the source term $q(x)$ from the measurement of $y$ on the time interval $(0,T)$ if $T$ is large enough.  
Our aim here is to propose well-chosen observers which allow, using a back-and-forth procedure, to retrieve $q$ from $y$ in the minimal observation time. More precisely, we prove the following result

\begin{claim}
\label{th1}
We can construct efficient observers for which the back-and-forth algorithm is convergent and which allow, using the measurement output $y(t)$ over the time-interval $(0,2)$, to reconstruct the unknown source term $q(x)$.
\end{claim}
\begin{remark}
Let us point out that since the spatial domain is given by $(0,1)$, the minimum observability time is given by $T=2$.
\end{remark}

 Note that, whenever the whole initial state $(u(0,.), u_t(0,.),q)$ is unknown, system \eqref{eq1} is not observable. In order to realize this, one can consider the simpler case where the source term $q$ is given by only the two first modes of the wave equation, $q(x)=q_1\sin(\pi x)+q_2\sin(2\pi x)$, where $q_1$ and $q_2$ are the unknown scalars. In this case, \eqref{eq1} becomes equivalent to two independent oscillators with different frequencies and with two unknown source terms $q_1\cos(\omega t)$ and $q_2\cos(\omega t)$. Moreover, the output is given by a linear combination of the position of the oscillators. This six dimensional system with one output is not observable. However, if we know the initial state of the oscillators, the two parameters $q_1$ and $q_2$ become identifiable.\\
The back-and-forth estimator allows us to take into account this knowledge of the initial state $(u(0,.), u_t(0,.))=(0,0)$ of the wave equation. On the contrary, if we had only used a forward observer, we would have lost the information on the initial state of the system and therefore there would have been no reason for the observer to converge to the real parameters.\\
 We prove the convergence of the algorithm using Lyapunov techniques and LaSalle's principle. One of the main difficulties comes from the fact that the precompactness of the trajectories is not ensured since we deal with an infinite dimensional system and we also have to use sharp mathematical estimates.\\

The paper is organized as follows.\\
In Section 2, we prove the equivalence between \eqref{eq1} and another system which consists in a system composed of a wave equation without any source term and an oscillator in cascade.  Here the unknown term to retrieve is the initial condition.
Section 3 is devoted to the proof of Claim \ref{th1} and is divided in two parts. In the first one we study the well-posedness of the proposed observer and in the second one we prove its convergence through the back-and-forth procedure. 
In Section 4, we provide some numerical simulations illustrating the efficiency of the method. 
Finally in Section 5, we propose and analyze the extension of the method to the case of wave equations with $N$-dimensional spatial domain. 

%--------------------------------------------------------------------------------
\section{An equivalent estimation problem}
We begin with  introducing the following system,

\begin{equation}
\label{eq2}
\left\{
\begin{array}{l}
w_{tt}-w_{xx}=0,\,(t,x)\in(0,T)\times(0,1),\\
w(t,0)=w(t,1)=0,\,t\in(0,T),\\
w(0,x)=q(x),\,w_t(0,x)=0,\,x\in(0,1),\\
\end{array}
\right. 
\end{equation}
\begin{equation}
\label{eq2osc}
\left\{
\begin{array}{l}
\dot{z_1}(t)=z_2(t),\,t\in(0,T),\\
\dot{z_2}(t)=-\omega^2z_1(t)+w_x(t,0),\,t\in(0,T),\\
z_1(0)=y(0),\, z_2(0)=\dot{y}(0),\\
Y(t)=z_1(t),\,t\in(0,T),
\end{array}
\right. 
\end{equation}
where $(w,w_t,z_1,z_2)$ represents the state of the system and $Y$ is the output. Here, the fact that $y(0)$ and $\dot{y}(0)$ are well-defined is a consequence of the definition of $y$ (see \eqref{eq1}) and the regularity property \eqref{eq00}. System \eqref{eq2}-\eqref{eq2osc} is nothing but a homogeneous wave equation and an oscillator in cascade. The unknown $q$ is now the initial datum of the system.
The aim of this section is to prove both following results
\begin{proposition}
\label{prop0}
There exists a unique solution $$(w,w_t,z_1,z_2)\in C([0,T];H^1_0(0,1))\times C([0,T];L^2(0,1))\times H^2(0,T)\times H^1(0,T)$$ to \eqref{eq2}. Moreover, $w$ satisfies the following hidden property
\begin{gather}
\label{extra}
w_x(t,0)\in L^2(0,T).
\end{gather}
\end{proposition}
and 
\begin{proposition}
\label{prop1}
For any $T>0$, for any $q\in H^2(0,1)$, $$y= Y \text{ in } H^2(0,T),$$ where $y$ denotes the output of system \eqref{eq1} and $Y$ denotes the output of system \eqref{eq2}.
\end{proposition}

We start by proving the Proposition~\ref{prop0} which ensures the well-posedness of the system \eqref{eq2}-\eqref{eq2osc}.

%----------------------------------------------------------------------------------------------------------------
\textbf{Proof of Proposition \ref{prop0}.} 
It follows from \cite[Lemma 3.6\, page 39]{Lionsbook88}  that there exists a unique $w\in C([0,T];H^1_0(0,1))\cap C^1([0,T];L^2(0,1))$ solution of 
\begin{equation}
\label{eq0f}
\left\{
\begin{array}{l}
w_{tt}-w_{xx}=0,\,(t,x)\in(0,T)\times(0,1),\\
w(t,0)=w(t,1)=0,\,t\in(0,T),\\
w(0,x)=q(x),\,w_t(0,x)=0,\,x\in(0,1).
\end{array}
\right. 
\end{equation}
From the hidden regularity property (see \cite[Th\'eor\`eme 4.1\, page 44]{Lionsbook88}),
\begin{gather}
\label{eq7}
w_x(t,0)\in L^2(0,T).
\end{gather}
Defining the matrix $A:=\left(\begin{array}{cc}0 & 1 \\-\omega^2 & 0\end{array}\right)$, the unique solution of the equation~\eqref{eq2osc} is given by, 
\begin{gather}
\label{eq8}
\left(\begin{array}{ccc}z_1(t) \\z_2(t)\end{array}\right)=e^{tA}\left(\begin{array}{c}y(0)\\\dot{y}(0)\end{array}\right)+\int_0^te^{(t-s)A}\left(\begin{array}{c}0 \\w_x(t,0)\end{array}\right).
\end{gather}
Using \eqref{eq8}, one easily gets
\begin{gather*}
\label{eq9}
\vert (z_1,z_2)\vert_{L^2(0,T)^2}\le C(1+\vert w_x(t,0)\vert_{L^2(0,T)}),
\end{gather*}
where $C$ denotes a positive constant which only depends on $\omega$ and $T$. 
Using first and  second equations in \eqref{eq2osc}, together with \eqref{eq7} and the last inequality we get 
\begin{gather}
\label{eq0g}
z_1  \text{ and }  z_2\in H^1(0,T),
\end{gather}
which implies, from the first equation in \eqref{eq2osc},  $z_1\in H^2(0,T)$. This ends the proof of Proposition \ref{prop0}.\cqfd
Now we start with the proof of the proposition~\ref{prop1} which ensures the equivalence (in the sense of the measurement output) between the equation~\eqref{eq1} and the system~\eqref{eq2}-\eqref{eq2osc}. 
%----------------------------------------------------------------------------------------------------------------

\textbf{Proof of Proposition \ref{prop1}.} 
We first prove the following lemma
\begin{lemma}
\label{lem0}
The solution $u$ to \eqref{eq1} satisfies
\begin{gather}
\label{eq0}
u\in C([0,T];H^3(0,1)\cap H^1_0(0,1)),\,u_t\in C([0,T];H^2(0,1)\cap H^1_0(0,1)), 
\end{gather}
and
\begin{gather}
\label{eq00}
u_x(t,0)\in H^2(0,T).
\end{gather}
\end{lemma}
\textbf{Proof of Lemma \ref{lem0}.} 
We proceed as in \cite[Appendix \, page 495]{1995-Yamamoto}. We recall that from \cite {1995-Yamamoto} , we already know that $u\in C^1([0,T];H^1_0(0,1))\cap C^2([0,T];L^2(0,1))$ and $u_x(t,0)\in H^1(0,T)$ (see in particular \cite[(1.3) page 183]{1995-Yamamoto}). We consider now the following system, 
 \begin{equation}
\label{eq0a}
\left\{
\begin{array}{l}
v_{tt}-v_{xx}=-\omega^2q(x)\cos(\omega t),\,(t,x)\in(0,T)\times(0,1),\\
v(t,0)=v(t,1)=0,\,t\in(0,T),\\
v(0,x)=q(x),\,v_t(0,x)=0,\,x\in(0,1).
\end{array}
\right. 
\end{equation}
Since $q\in H^2(0,1)\cap H^1_0(0,1)$ it follows from \cite[Lemma 3.6\, page 39]{Lionsbook88} that 
\begin{gather}
\label{eq0b}
v\in C([0,T];H^1_0(0,1))\cap C^1([0,T];L^2(0,1)),
\end{gather}
and 
\begin{gather}
\label{eq0bb}
v_x(t,0)\in L^2(0,T).
\end{gather}
We can set $$\tilde{u}(t,x):=\int_0^t\int_0^sv(\tau,x)\,d\tau\,ds,\, t\in(0,T),\,x\in(0,1),$$ so that $\tilde{u}\in C^2([0,T];H^1_0(0,1))\cap C^3([0,T];L^2(0,1))$ and $\tilde{u}$ satisfies \eqref{eq1}. By uniqueness of weak solution to \eqref{eq1}, we get $$u(t,x)=\tilde{u}(t,x) \text{ and } u_{tt}(t,x)=v(t,x),\, t\in(0,T),\,x\in (0,1).$$ Consequently 
\begin{gather}
\label{eq0c}
u\in C^2([0,T];H^1_0(0,1))\cap C^3([0,T];L^2(0,1)),
\end{gather}
and from \eqref{eq0bb},
\begin{gather*}
\label{eq0cc}
u_x(t,0)\in H^2(0,T).
\end{gather*}
 Now, using the fact that $u$ satisfies the equation~\eqref{eq1} and that $u_{tt}\in C([0,T];H^1_0(0,1))$ (see \eqref{eq0c}) and that $q\in H^2(0,1)\cap H^1_0(0,1)$, we obtain 
 \begin{gather*}
\label{eq0d}
u_{xx}\in C([0,T];H^1_0(0,1)),
\end{gather*}
and thus 
 \begin{gather*}
\label{eq0e}
u\in C([0,T];H^3\cap H^1_0(0,1)).
\end{gather*}
Using one more time the fact that $u$ satisfies~\eqref{eq1}, that  $q\in H^2(0,1)\cap H^1_0(0,1)$ and that $u_{tt}\in C^1([0,T];L^2(0,1))$ (see \eqref{eq0c}), it follows 
\begin{gather}
\label{eq0eb}
u_{xx}\in C^1([0,T];L^2(0,1)).
\end{gather}
 Consequently, it follows from \eqref{eq0c} and \eqref{eq0eb} that $$u\in C^1([0,T];H^2(0,1)\cap H^1_0(0,1)).$$ This ends the proof of Lemma \ref{lem0}.\cqfd
 % - -- - - - --- - -- - -  - -- - -  - -- - - - - - - - - - - - - - - - - - - - - - - 
Since $u\in C([0,T];H^3(0,1)\cap H^1_0(0,1))\cap C^1([0,T];H^2(0,1)\cap H^1_0(0,1))$ and $q\in H^2(0,1)\cap H^1_0(0,1)$, we can define $W\in C([0,T],H^1_0(0,1))\cap C^1([0,T];L^2(0,1))$ as 
\begin{gather}
\label{defW}
W:=\omega^2u+u_{tt}=\omega^2u+u_{xx}+\cos(\omega t)q(x).
\end{gather}
 It is easy to see that $W$ satisfies 
\begin{equation}
\label{eq0i}
\left\{
\begin{array}{l}
W_{tt}-W_{xx}=0,\,(t,x)\in(0,T)\times(0,1),\\
W(t,0)=W(t,1)=0,\,t\in(0,T),\\
W(0,x)=q(x),\,W_t(0,x)=0,\,x\in(0,1),
\end{array}
\right. 
\end{equation}
and thus, by uniqueness of the weak solution of the previous system, $$w=W \text{ in }  C([0,T],H^1_0(0,1))\cap C^1([0,T];L^2(0,1))$$ and $$w_x(t,0)=W_x(t,0)\text{ in } L^2(0,T).$$ In other words, 
\begin{gather}
\label{eq11}
w=\omega^2u+u_{xx}+\cos(\omega t)q(x) \text{ in } C([0,T],H^1_0(0,1))\cap C^1([0,T];L^2(0,1)),
\end{gather}
and one can compute from \eqref{defW}  (computations making sense from \eqref{eq00}),
\begin{gather*}
w_x(t,0)=\omega^2u_x(t,0)+u_{ttx}(t,0) \text{ in } L^2(0,T),
\end{gather*}
and thus from \eqref{eq1}
\begin{gather}
\label{tchoum}
w_x(t,0)= \omega^2y(t)+\ddot{y}(t)\text{ in } L^2(0,T).
\end{gather}
Now, from \eqref{eq2osc}, 
\begin{gather}
\label{was}
\ddot{Y}+\omega^2 Y=w_x(t,0)\text{ in } L^2(0,T).
\end{gather}
Finally, from \eqref{eq2osc}, \eqref{tchoum}, \eqref{was},
\begin{equation}
\label{eq12}
\left\{
\begin{array}{l}
\ddot{Y}+\omega^2Y=\ddot{y}+\omega^2y \text{ in } L^2(0,T),\\
Y(0)=y(0),\,\dot{Y}(0)=\dot{y}(0),
\end{array}
\right. 
\end{equation}
which concludes the proof of Proposition \ref{prop1} by uniqueness of the solution to the ordinary differential equation $\ddot{\xi}+\omega^2\xi=0,\,\xi(0)=\dot{\xi}(0)=0.$ 
\cqfd
\begin{remark}
\label{remrem}
One can note that using usual energy estimates and  the fact that  $q\in H^2(0,1)\cap H^1_0(0,1)$, one can prove more regularity for the solution $(w,w_t)$ to \eqref{eq2}. Precisely, one actually obtains $$w\in C([0,T],H^2(0,1)\cap H^1_0(0,1))\cap C^1([0,T];H^1_0(0,1)).$$ 
\end{remark}

As a consequence of Proposition \ref{prop1}, we will now focus on system \eqref{eq2}-\eqref{eq2osc} and provide well-chosen observers for this sytem in order to prove Claim \ref{th1}.
%--------------------------------------------------------------------------------------------------------------
\section{Observer design}
As we have already mentioned in the introduction, we are going to prove Claim \ref{th1} using back-and-forth observers. We first rewrite system \eqref{eq2}-\eqref{eq2osc} in the following way
\begin{equation}
\label{eq2*}
\left\{
\begin{array}{l}
w^1_{t}=w^2,\,(t,x)\in(0,T)\times(0,1),\\
w^2_{t}=w^1_{xx},\,(t,x)\in(0,T)\times(0,1),\\
w^1(t,0)=w^1(t,1)=0,\,t\in(0,T),\\
w^1(0,x)=q(x),\,w^2(0,x)=0,\,x\in(0,1),\\
\end{array}
\right. 
\end{equation}
\begin{equation}
\label{eq2osc*}
\left\{
\begin{array}{l}
\dot{z_1}(t)=z_2(t),\,t\in(0,T),\\
\dot{z_2}(t)=-\omega^2z_1(t)+w^1_x(t,0),\,t\in(0,T),\\
z_1(0)=y(0),\, z_2(0)=\dot{y}(0),\\
Y(t)=z_1(t),\,t\in(0,T),
\end{array}
\right. 
\end{equation}
where, from above results, $$(w^1,w^2,z^1,z^2)\in C([0,T];H^2(0,1)\cap H^1_0(0,1))\times C([0,T];H^1_0(0,1))\times H^2(0,1)\times H^1(0,1).$$
Let us recall that our aim is to propose an algorithm which allows to retrieve the unknown $q$ from the measurement output $Y$ on the time interval $(0,T)$. Our idea consists in designing a well-chosen asymptotic observer, which ensures the decrease of a same Lyapunov function in the back-and-forth procedure.  To this aim, we first make the previous system periodical. More precisely, we define $(W^1,W^2,Z^1,Z^2)$ on $((0,+\infty)\times(0,1))^2\times (0,T)^2$ as the solution  to the following periodical system
\begin{equation}
\label{eq2**}
\left\{
\begin{array}{l}
W^1_{t}=W^2,\,(t,x)\in(2kT,(2k+1)T)\times(0,1),\vspace{0,1cm}\\
W^2_{t}=W^1_{xx},\,(t,x)\in(2kT,(2k+1)T)\times(0,1),\vspace{0,1cm}\\
W^1_{t}=-W^2,\,(t,x)\in((2k+1)T,(2k+2)T)\times(0,1),\vspace{0,1cm}\\
W^2_{t}=-W^1_{xx},\,(t,x)\in((2k+1)T,(2k+2)T)\times(0,1),\vspace{0,1cm}\\
W^1(t,0)=W^1(t,1)=0,\,t\in(0,+\infty),\vspace{0,1cm}\\
W^1(0,x)=q(x),\,W^2(0,x)=0,\,x\in(0,1),
\end{array}
\right. 
\end{equation}
\begin{equation}
\label{eq2osc**}
\left\{
\begin{array}{l}
\dot{Z_1}(t)=Z_2(t),\,t\in(2kT,(2k+1)T),\vspace{0,1cm}\\
\dot{Z_2}(t)=-\omega^2Z_1(t)+W^1_x(t,0),\,t\in(2kT,(2k+1)T),\vspace{0,1cm}\\
\dot{Z_1}(t)=-Z_2(t),\,t\in((2k+1)T,(2k+2)T),\vspace{0,1cm}\\
\dot{Z_2}(t)=\omega^2Z_1(t)-W^1_x(t,0),\,t\in((2k+1)T,(2k+2)T),\vspace{0,1cm}\\
Z_1(0)=y(0),\, Z_2(0)=\dot{y}(0),
\end{array}
\right. 
\end{equation}
for  $k\in \mathbb{N}$.
Using previous results, one easily sees that $$(W^1,W^2)\in L^{\infty}_{loc}(\mathbb{R}_+;H^2(0,1)\cap H^1_0(0,1))\times L^{\infty}_{loc}(\mathbb{R}_+;H^1_0(0,1)),$$ 
and $$(Z^1,Z^2)\in H^2_{loc}(\mathbb{R}_+)\times H^1_{loc}(\mathbb{R}_+).$$
Indeed, the above periodic system represents the system \eqref{eq2*}-\eqref{eq2osc*} on the time intervals $(2kT,(2k+1)T)$, and the same system in the time-reversed manner on the intervals $((2k+1)T,(2k+2)T)$. In particular, one has
\begin{equation*}
\left\{
\begin{array}{l}
Y(t)=Y(t-2kT),\,t\in(2kT,(2k+1)T),\vspace{0,1cm}\\
Y(t)=Y((2k+2)T-t),\,t\in((2k+1)T,(2k+2)T).
\end{array}
\right. 
\end{equation*}

%%%%%%%%%%%%%%%%%%%%%%%
The asymptotic observer we propose is the following

\begin{equation}
\label{eq13}
\left\{
\begin{array}{l}
\hat{W}^1_{t}=\hat{W}^2,\,(t,x)\in(2kT,(2k+1)T)\times(0,1),\vspace{0,1cm}\\
\hat{W}^2_{t}=\hat{W}^1_{xx},\,(t,x)\in(2kT,(2k+1)T)\times(0,1),\vspace{0,1cm}\\
\hat{W}^1_{t}=-\hat{W}^2,\,(t,x)\in((2k+1)T,(2k+2)T)\times(0,1),\vspace{0,1cm}\\
\hat{W}^2_{t}=-\hat{W}^1_{xx},\,(t,x)\in((2k+1)T,(2k+2)T)\times(0,1),\vspace{0,1cm}\\
\hat{W}^{1}(t,0)=\gamma_1(\hat{Z}_1(t)-Y(t))+\gamma_1\gamma_2(\hat{Z}_3(t)-\int_0^{t}Y(s)\,ds),\,t\in(0,\infty),\vspace{0,1cm}\\
\hat{W}^1(t,1)=0,\,t\in(0,+\infty),\vspace{0,1cm}\\
\hat{W}^1(0,x)=\hat{W}^2(0,x)=0,\,x\in(0,1),
\end{array}
\right. 
\end{equation}
\begin{equation}
\label{eq14}
\left\{
\begin{array}{l}
\dot{\hat{Z}}_1(t)=\hat{Z}_2(t)-\gamma_2(\hat{Z}_1(t)-Y(t)),\,t\in(2kT,(2k+1)T),\vspace{0,1cm}\\
\dot{\hat{Z}}_2(t)=-\omega^2\hat{Z}_1(t)+\hat{W}^1_x(t,0),\,\phantom{ggg}t\in(2kT,(2k+1)T),\vspace{0,1cm}\\
\dot{\hat{Z}}_1(t)=-\hat{Z}_2(t)-\gamma_2(\hat{Z}_1(t)-Y(t)),\,t\in((2k+1)T,(2k+2)T),\vspace{0,1cm}\\
\dot{\hat{Z}}_2(t)=\omega^2\hat{Z}_1(t)-\hat{W}^1_x(t,0),\,\phantom{ttffft}t\in((2k+1)T,(2k+2)T),\vspace{0,1cm}\\
\dot{\hat{Z}}_3(t)=\hat{Z}_1(t),\,t\in(0,\infty),\vspace{0,1cm}\\
\hat{Z}_1(0)=\hat{Z}_2(0)=\hat{Z}_3(0)=0.
\end{array}
\right. 
\end{equation}
Here $\gamma_1$ and $\gamma_2$, the observer gains, are strictly positive constants to be fixed as the design parameters.
Before studying the well-posedness of system \eqref{eq13}-\eqref{eq14}, let us introduce error  equations. The error term being defined as the difference between the observer and the observed system, error  equations are the following
\begin{equation}
\label{eq17}
\left\{
\begin{array}{l}
\tilde{W}^1_{t}=\tilde{W}^2,\,(t,x)\in(2kT,(2k+1)T)\times(0,1),\vspace{0,1cm}\\
\tilde{W}^2_{t}=\tilde{W}^1_{xx},\,(t,x)\in(2kT,(2k+1)T)\times(0,1),\vspace{0,1cm}\\
\tilde{W}^1_{t}=-\tilde{W}^2,\,(t,x)\in((2k+1)T,(2k+2)T)\times(0,1),\vspace{0,1cm}\\
\tilde{W}^2_{t}=-\tilde{W}^1_{xx},\,(t,x)\in((2k+1)T,(2k+2)T)\times(0,1),\vspace{0,1cm}\\
\tilde{W}^{1}(t,0)=\gamma_1\tilde{Z}_1(t)+\gamma_1\gamma_2\tilde{Z}_3(t),\,t\in(0,\infty),\vspace{0,1cm}\\
\tilde{W}^1(t,1)=0,\,t\in(0,+\infty),\vspace{0,1cm}\\
\tilde{W}^1(0,x)=-q(x),\,\tilde{W}^2(0,x)=0,\,x\in(0,1),
\end{array}
\right. 
\end{equation}
\begin{equation}
\label{eq18}
\left\{
\begin{array}{l}
\dot{\tilde{Z}}_1(t)=\tilde{Z}_2(t)-\gamma_2\tilde{Z}_1(t),\,t\in(2kT,(2k+1)T),\vspace{0,1cm}\\
\dot{\tilde{Z}}_2(t)=-\omega^2\tilde{Z}_1(t)+\tilde{W}^1_x(t,0),\,t\in(2kT,(2k+1)T),\vspace{0,1cm}\\
\dot{\tilde{Z}}_1(t)=-\tilde{Z}_2(t)-\gamma_2\tilde{Z}_1(t),\,t\in((2k+1)T,(2k+2)T),\vspace{0,1cm}\\
\dot{\tilde{Z}}_2(t)=\omega^2\tilde{Z}_1(t)-\tilde{W}^1_x(t,0),\,t\in((2k+1)T,(2k+2)T),\vspace{0,1cm}\\
\dot{\tilde{Z}}_3(t)=\tilde{Z}_1(t),\,t\in(0,\infty),\vspace{0,1cm}\\
\tilde{Z}_1(0)=-y(0),\,\tilde{Z}_2(0)=-\dot{y}(0),\,\tilde{Z}_3(0)=0.
\end{array}
\right. 
\end{equation}

%--------------------------------------------------------------------------------------
\subsection{Well-posedness}
From now on, and till the end, we denote $H^1_r(0,1):=\{v\in H^1(0,1)\text{ s.t. } v(1)=0\}.$ Let us assume for the moment that the following proposition holds (for the proof we refer the reader to the Appendix). 
\begin{proposition}
\label{prop2}
For any $T>0$, for any $\gamma_1$ and $\gamma_2>0$, for any $(q^0,q^1)\in H^2(0,1)\cap H^1_0(0,1)\times H^1_0(0,1)$ and any $(\xi_1^0,\xi_2^0,\xi_3^0)\in\mathbb{R}^3$, there exists a unique solution $$(v^1,v^2)\in C([0,+\infty);H^2(0,1)\cap H^1_r(0,1))\times C([0,+\infty);H^1_r(0,1)),$$
$$(\xi_1,\xi_2,\xi_3)\in H^1([0,+\infty))\times L^2([0,+\infty))\times H^2([0,+\infty))$$

 to the following periodical Cauchy problem
\begin{equation}
\label{eqerror}
\left\{
\begin{array}{l}
v^1_{t}=v^2,\,(t,x)\in(2kT,(2k+1)T)\times(0,1),\vspace{0,1cm}\\
v^2_{t}=v^1_{xx},\,(t,x)\in(2kT,(2k+1)T)\times(0,1),\vspace{0,1cm}\\
v^1_{t}=-v^2,\,(t,x)\in((2k+1)T,(2k+2)T)\times(0,1),\vspace{0,1cm}\\
v^2_{t}=-v^1_{xx},\,(t,x)\in((2k+1)T,(2k+2)T)\times(0,1),\vspace{0,1cm}\\
v^1(t,1)=0,\,v^1(t,0)=\gamma_1\xi_1(t)+\gamma_1\gamma_2\xi_3(t),\,t\in((2k+1)T,(2k+2)T),\vspace{0,1cm}\\
v^1(0,x)=q^0(x),\,v^2(0,x)=q^1(x),\,x\in(0,1),
\end{array}
\right. 
\end{equation}
\begin{equation}
\label{eqerrorbis}
\left\{
\begin{array}{l}
\dot{\xi_1}(t)=\xi_2(t)-\gamma_2\xi_1(t),\,t\in(2kT,(2k+1)T),\vspace{0,1cm}\\
\dot{\xi_2}(t)=-\omega^2\xi_1(t)+v^1_x(t,0),\,t\in(2kT,(2k+1)T),\vspace{0,1cm}\\
\dot{\xi_1}(t)=-\xi_2(t)-\gamma_2\xi_1(t),\,t\in((2k+1)T,(2k+2)T),\vspace{0,1cm}\\
\dot{\xi_2}(t)=\omega^2\xi_1(t)-v^1_x(t,0),\,t\in((2k+1)T,(2k+2)T),\vspace{0,1cm}\\
\dot{\xi_3}(t)=\xi_1(t),\,t\in(0,\infty),\vspace{0,1cm}\\
\xi_1(0)=\xi_1^0,\, \xi_2(0)=\xi_2^0,\, \xi_2(0)=\xi_3^0,
\end{array}
\right. 
\end{equation}
where $k\in \mathbb{N}$.
Moreover $(v^1,v^2,\xi_1,\xi_2,\xi_3)$ satisfy the following energy identities: for any $t\in (0,+\infty)$,
\begin{eqnarray}
\begin{array}{l}
\label{eq26}
\vert v^1_x(t,.)\vert^2_{L^2(0,1)}+\vert v^2(t,.)\vert^2_{L^2(0,1)}+\gamma_1\vert\xi_2(t)\vert^2+\gamma_1\omega^2\vert\xi_1(t)\vert^2+2\gamma_1\gamma_2\omega^2\vert\xi_1\vert^2_{L^2(0,t)}=\vspace{0,2cm}\\
\phantom{ttttttdkkmmmmmmmkkkkkkkkkkkddt}\vert q^0_x\vert^2_{L^2(0,1)}+\vert q^1\vert^2_{L^2(0,1)}+\gamma_1\vert\xi_2^0\vert^2+\gamma_1\omega^2\vert\xi_1^0\vert^2.\end{array}
\end{eqnarray}
and 
\begin{eqnarray}
\begin{array}{l}
\label{eq26bis}
\displaystyle \frac{1}{2}(\vert v^2_{t}(t,.)\vert^2_{L^2(0,1)}+\vert v^2_{x}(t,.)\vert^2_{L^2(0,1)}+\gamma_1\omega^4\vert\xi_1(t)\vert^2)+\frac{1}{4}\gamma_1\vert v^1_x(t,0)\vert^2
\vspace{0,1cm}\\
\displaystyle\phantom{ttttt}+\frac{1}{2}\gamma_1\gamma_2\omega^2\int_0^t\vert \xi_2(s)\vert^2\,ds+\gamma_1\omega^2\vert \xi_2(t)\vert^2\le C(\vert q^0\vert^2_{H^2(0,1)}+\vert q^1\vert^2_{H^1(0,1)}+\vert \xi_1^0\vert^2+\vert \xi_2^0\vert^2),
\end{array}
\end{eqnarray}
where $C_{\omega}$ denotes a positive constant which only depends on $\omega$.
% there exists $C>0$ which only depends on $\omega$ and does not depend on the time $T$ such that
\end{proposition}

%\begin{remark}
%We recall that for any $f\in H^1(0,T)$, the solution $(v^1,v^2)$ to 
%\begin{equation}
%\label{eq27}
%\left\{
%\begin{array}{l}
%v^1_{t}=v^2,\,(t,x)\in(0,T)\times(0,1),\\
%v(t,1)=0,\,v(t,0)=f(t),\,t\in(0,T),\\
%v(0,x)=v^0(x),\,v_t(0,x)=v^1(x),\,x\in(0,1),
%\end{array}
%\right. 
%\end{equation}
%can be defined by transposition  (see \cite[page 47]{Lionsbook88}).
%\end{remark}
One easily sees that the well-posedness of system \eqref{eq17}-\eqref{eq18} in $$C([0,+\infty);H^2(0,1)\cap H^1_r(0,1))\times C([0,+\infty);H^1_r(0,1))\times H^1([0,+\infty))\times L^2([0,+\infty))\times H^2([0,+\infty))$$
directly follows from Proposition \ref{prop2}. The well-posedness of system \eqref{eq13}-\eqref{eq14} in $$L^{\infty}_{loc}(\mathbb{R}_+;H^2(0,1)\cap H^1_r(0,1))\times L^{\infty}_{loc}(\mathbb{R}_+;H^1_r(0,1))\times H^1_{loc}(\mathbb{R}_+)\times L^2_{loc}(\mathbb{R}_+)\times H^2_{loc}(\mathbb{R}_+)$$ then follows (using also the well-posedness of system \eqref{eq2**}-\eqref{eq2osc**} in the same space).

%-----------------------------------------------------------------------------------------------------------
%-----------------------------------------------------------------------------------------------------------

\subsection{Asymptotic analysis}
The main goal of this section is to prove the following result, implying Claim \ref{th1}
\begin{theorem}
\label{th2}
 For any $T\ge 2$,
\begin{gather}
\label{th2eq1}
\lim\limits_{n\rightarrow +\infty}(\tilde{W}^1(2nT,.),\tilde{W}^2(2nT,.))= (0,0)\text{ in } H^1_r(0,1)\times L^2(0,1),\\
\label{th2eq2}
\lim\limits_{n\rightarrow +\infty}(Z_1(2nT),Z_2(2nT),Z_3(2nT))= (0,0,0).
\end{gather}
\end{theorem}
\textbf{Proof of Theorem \ref{th2}.}Let us assume that \eqref{th2eq1} and \eqref{th2eq2} do not hold. Then, there exist a positive constant $\alpha$ and  a subsequence $(\phi(n))_{n\ge0}\in \mathbb{N}^\mathbb{N}$  with $\lim\limits_{n\rightarrow +\infty}\phi(n)=+\infty$ such that for any $n\ge0$,
\begin{equation}
\label{th2eq3}
\max\left(\|\tilde{W}^1(t,.)\|_{H^1_r(0,1)},\|\tilde{W}^2(t,.)\|_{L^2(0,1)}, |\tilde Z^1(t)|,|\tilde Z^2(t)|,|\tilde Z^3(t)|\right)\Big|_{t=2\phi(n)T}> \alpha.
\end{equation}\normalsize

From Proposition \ref{prop2}, $$(\tilde{W}^1,\tilde{W}^2)\in L^{\infty}(\mathbb{R}_+;H^2(0,1)\cap H^1_r(0,1))\times L^{\infty}(\mathbb{R}_+;H^1_r(0,1)),$$ 
$$(\tilde{Z}_1,\tilde{Z}_2,\tilde{Z}_3)\in H^1(\mathbb{R}_+)\times L^2(\mathbb{R}_+)\times H^2(\mathbb{R}_+),$$

and for any time $t\in \mathbb{R}_+^*$, 
\begin{eqnarray}
\begin{array}{l}
\label{eq70}
\vert \tilde{W}^1_x(t,.)\vert^2_{L^2(0,1)}+\vert \tilde{W}^2(t,.)\vert^2_{L^2(0,1)}+\gamma_1\vert\tilde{Z}_2(t)\vert^2
+\gamma_1\omega^2\vert\tilde{Z}_1(t)\vert^2+2\gamma_1\gamma_2\omega^2\vert\tilde{Z}_1\vert^2_{L^2(0,t)}=\vspace{0,2cm}\\
\phantom{tttttllllmmmmmmmmmmmmmmtttttttttttlt}\vert q\vert^2_{H^1(0,1)}+\gamma_1\vert \dot{y}(0)\vert^2+\gamma_1\omega^2\vert y(0)\vert^2.
\end{array}
\end{eqnarray}
Since $\{(\tilde{W}^1(2\phi(n)T,.),\tilde{W}^2(2\phi(n)T,.)),\,n\ge0\}$ is bounded in $H^2(0,1)\times H^1_r(0,1)$, it follows from Rellich's theorem that there exists a subsequence of $(\phi(n))_{n\ge0}$, that, for convenience, we still denote $(\phi(n))_{n\ge0}$, and there exists $(W^{\infty,1},W^{\infty,2})\in H^1_r(0,1)\times L^2(0,1)$ such that 
\begin{gather}
\label{cv}
\lim\limits_{n\rightarrow +\infty}(\tilde{W}^1(2\phi(n)T,.),\tilde{W}^2(2\phi(n)T,.))=(W^{\infty,1},W^{\infty,2}) \text{ in } H^1_r(0,1)\times L^2(0,1).
\end{gather}
From \eqref{eq70},  $\tilde{W}^1(2\phi(n)T,0)$, $\tilde{Z}_1(2\phi(n)T)$ and $\tilde{Z}_2(2\phi(n)T)$ are bounded, and from third and sixth equations in \eqref{eq17} so is $\tilde{Z}_3(2\phi(n)T)$. Thus, (up to a subsequence of $(\phi(n))_{n\ge0}$),  there exists $(Z_1^{\infty},Z_2^{\infty},Z_3^{\infty})\in \mathbb{R}^3$ such that
\begin{gather}
\label{cvbis}
\lim\limits_{n\rightarrow +\infty} (\tilde{Z}_1(2\phi(n)T),\tilde{Z}_2(2\phi(n)T),\tilde{Z}_3(2\phi(n)T))=(Z_1^{\infty},Z_2^{\infty},Z_3^{\infty}).
\end{gather}

Let now $$(v^1,v^2,x_1,x_2,x_3)\in C([0,T];H^1_r(0,1))\times C([0,T];L^2(0,1))\times L^2(0,T)\times L^2(0,T)\times L^2(0,T)$$ be solution to
\begin{equation}
\label{eq80}
\left\{
\begin{array}{l}
v^1_{t}=v^2,\,(t,x)\in(0,T)\times(0,1),\\
v^2_{t}=v^1_{xx},\,(t,x)\in(0,T)\times(0,1),\\
v^1(t,1)=0,\,v^1(t,0)=\gamma_1x_1(t)+\gamma_1\gamma_2x_3(t),\,t\in(0,T),\\
\dot{x_1}(t)=x_2(t)-\gamma_2x_1(t),\,t\in(0,T),\\
\dot{x}_2(t)=-\omega^2x_1(t)+v^1_x(t,0),\,t\in(0,T),\\
\dot{x}_3(t)=x_1(t),\,t\in(0,T),\\
v^1(0,x)=W^{\infty,1}(x),\,v^2(0,x)=W^{\infty,2}(x),\,x\in(0,1),\\
x_1(0)=Z_1^{\infty},\,x_2(0)=Z_2^{\infty},\,x_3(0)=Z_3^{\infty}.
\end{array}
\right. 
\end{equation}
Finally, let us define, for any $n\ge0$, $$(v^{1,n},v^{2,n},x_1^n,x_2^n,x_3^n)\in C([0,T];H^1_r(0,1))\times C([0,T];L^2(0,1))\times L^2(0,T)\times L^2(0,T)\times L^2(0,T)$$ by
\begin{gather}
 \label{eq86}
(v^{1,n}(t,.),v^{2,n}(t,.)):=(\tilde{W}^1(2\phi(n)T+t,.),\tilde{W}^2(2\phi(n)T+t,.))
\end{gather}
and 
\begin{gather}
 \label{eq87}
(x_1^n(t),x_2^n(t),x_3^n(t)):=(\tilde{Z}_1(2\phi(n)T+t),\tilde{Z}_2(2\phi(n)T+t),\tilde{Z}_3(2\phi(n)T+t)),
\end{gather}
 for any $t\in(0,T)$. With such a definition, for any $n\ge0$, $(v^{1,n},v^{2,n},x_1^n,x_2^n,x_3^n)$ is solution to
\begin{equation}
\label{eq81}
\left\{
\begin{array}{l}
v^{1,n}_{t}=v^{2,n},\,(t,x)\in(0,T)\times(0,1),\vspace{0,1cm}\\
v^{2,n}_{t}=v^{1,n}_{xx},\,(t,x)\in(0,T)\times(0,1),\vspace{0,1cm}\\
v^{1,n}(t,1)=0,\,v^{1,n}(t,0)=\gamma_1x^n_1(t)+\gamma_1\gamma_2 x^n_3(t),\,t\in(0,T),\vspace{0,1cm}\\
\dot{x}_1^n(t)=x_2^n(t)-\gamma_2 x_1^n(t),\,t\in(0,T),\vspace{0,1cm}\\
\dot{x}^n_2(t)=-\omega^2x_1^n(t)+v^{1,n}_x(t,0),\,t\in(0,T),\vspace{0,1cm}\\
\dot{x}^n_3(t)=x_1^n(t),\,t\in(0,T),\vspace{0,1cm}\\
v^{1,n}(0,x)=\tilde{W}^1(2\phi(n)T,x),\,x\in(0,1),\\
v^{2,n}(0,x)=\tilde{W}^2(2\phi(n)T,x),\,x\in(0,1),\vspace{0,1cm}\\
x_1^n(0)=\tilde{Z}_1(2\phi(n)T),\,x_2^n(0)=\tilde{Z}_2(2\phi(n)T),\,x_3^n(0)=\tilde{Z}_3(\phi(n)T).\end{array}
\right. 
\end{equation}
From \eqref{eq26}, \eqref{cv}-\eqref{eq80} and \eqref{eq81} (more precisely by continuity of flow with respect to the initial state),
\begin{gather}
 \label{eq84}
\lim\limits_{n\rightarrow +\infty}(v^{1,n},v^{2,n})= (v^1,v^2)\text{ in } L^{\infty}((0,T);H^1_r(0,1))\times L^{\infty}((0,T);L^2(0,1)),\\
\label{eq84bis}
\lim\limits_{n\rightarrow +\infty}(x_1^n,x_2^n,x_3^n)=(x_1,x_2,x_3)\text{ in } L^{\infty}(0,T)^3.
\end{gather}

We now introduce the following Lyapunov function
\begin{eqnarray}
\begin{array}{l}
 \label{eq82}
\mathcal{V}(\tilde{W}^1,\tilde{W}^2,\tilde{Z}_1,\tilde{Z}_2,\tilde{Z}_3)(t):=\displaystyle\frac{1}{2}\Big(\int_0^1(\vert \tilde{W}^1_x(t,x)\vert^2+\vert \tilde{W}^2(t,x)\vert^2)\,dx\\
\phantom{nnnnnkkkkkkkkkkkkkkknnnn,,,,,,,,,,,,nn}+\gamma_1\omega^2\vert \tilde{Z}_1(t)\vert^2+\gamma_1\vert \tilde{Z}_2(t)\vert^2\Big),\,t\in(0,T).
\end{array}
\end{eqnarray}
Indeed, $\mathcal{V}(t)\ge0,\,t\ge0$ and using \eqref{eq17} and \eqref{eq18}, one can compute, for any time $t\in(0,T)$,
\begin{gather}
 \label{eq82bis}
\frac{d}{dt}\mathcal{V}(\tilde{W}^1,\tilde{W}^2,\tilde{Z}_1,\tilde{Z}_2,\tilde{Z}_3)(t)=-\gamma_1\gamma_2\omega^2\vert \tilde{Z}_1(t)\vert^2\le 0.
\end{gather}
The function $\mathcal{V}$ is positive, decreasing. Consequently there exists $l\ge0$ such that
\begin{gather}
\label{eq88}
\lim\limits_{t\rightarrow +\infty}\mathcal{V}(\tilde{W}^1,\tilde{W}^2,\tilde{Z}_1,\tilde{Z}_2,\tilde{Z}_3)(t)=l.
\end{gather}

On the other side,  from \eqref{eq86}, \eqref{eq87}, \eqref{eq84} and \eqref{eq82}, one has, for any $t\in(0,T)$,
\begin{eqnarray}
\begin{array}{l}
 \label{eq85}
  \lim\limits_{n\rightarrow+\infty}\mathcal{V}(\tilde{W}^1,\tilde{W}^2,\tilde{W}^2,\tilde{Z}_1,\tilde{Z}_2,\tilde{Z}_3)(2\phi(n)T+t)=\displaystyle\frac{1}{2}\Big(\int_0^1(\vert v^1_x(t,x)\vert^2\\
    \phantom{ttttttttttttttttttttttttttttt}+\vert v^2(t,x)\vert^2)\,dx
+\gamma_1\omega^2\vert x_1(t)\vert^2+\gamma_1\vert x_2(t)\vert^2\Big),
\end{array}
\end{eqnarray}

Thus it follows from \eqref{eq88} and \eqref{eq85} that for any time $t\in(0,T)$,
$$\mathcal{V}(v^1,v^2,x_1,x_2,x_3)(t)=l.$$
In other words, $t\mapsto\mathcal{V}(v^1,v^2,x_1,x_2,x_3)(t)$  is constant on $(0,T)$ and thus
\begin{gather}
 \label{eq89}
\frac{d}{dt}\mathcal{V}(v^1,v^2,x_1,x_2,x_3)(t)=-\gamma_1\gamma_2\omega^2\vert x_1(t)\vert^2=0,\,\qquad t\in(0,T).
\end{gather}
Therefore, we  obtain
\begin{gather}
 \label{eq90}
x_1\equiv0 \text{ on } (0,T).
\end{gather}
Consequently, from the fourth  and sixth lines in \eqref{eq80}, we also get
\begin{gather}
 \label{eq91}
x_2\equiv0,~\quad x_3\equiv  Z_3^\infty \quad \text{ on } (0,T),
\end{gather}
where $Z_3^\infty$ is a real constant.
Then, using the third and fifth lines in \eqref{eq80} and, \eqref{eq90} and \eqref{eq91} we see that
\begin{gather}
 \label{eq93}
v^1_x(t,0)\equiv0,\quad v^1(t,0)=\gamma_1\gamma_2 Z_3^\infty \quad\text{ on } (0,T).
\end{gather}
Consequently, \eqref{eq80} reduces to
\begin{equation}
\label{eq80bis}
\left\{
\begin{array}{l}
v^1_{t}=v^2,\,(t,x)\in(0,T)\times(0,1),\\
v^2_{t}=v^1_{xx},\,(t,x)\in(0,T)\times(0,1),\\
v^1(t,1)=v^1_x(t,0)=0,\,t\in(0,T),\\
v^1(t,0)=\gamma_1\gamma_2 Z_3^\infty ,\,t\in(0,T),\\
v^1(0,x)=W^{\infty,1}(x),\,v^2(0,x)=W^{\infty,2}(x),\,x\in(0,1),\\
x_1(t)=x_2(t)=0,~ x_3(t)=Z_3^\infty,\,t\in(0,T).
\end{array}
\right. 
\end{equation}

We recall that if $W^{\infty,1}(x)=\sum\limits_{k=1}^{+\infty}a_k\cos(\frac{(2k+1)\pi x}{2})$ and $W^{\infty,2}(x)=\sum\limits_{k=1}^{+\infty}b_k\cos(\frac{(2k+1)\pi x}{2})$, then $v^1$ can be expressed , for any $(t,x)\in(0,T)\times(0,1)$ as
\begin{multline*}
 v^1(t,x)=\sum_{k=1}^{\infty}\Big[a_k\cos\left(\frac{(2k+1)\pi t}{2}\right)+
 \frac{2b_k}{(2k+1)\pi }\sin\left(\frac{(2k+1)\pi t}{2}\right)\Big]\cos\left(\frac{(2k+1)\pi x}{2}\right).
\end{multline*}
At $x=0$, we have
\begin{multline*}
v^1(t,0)=\gamma_1\gamma_2 Z_3^{\infty}=\sum_{k=1}^{\infty}\Big[ a_k\cos\left(\frac{(2k+1)\pi t}{2}\right)
+\frac{2b_k}{(2k+1)\pi }\sin\left(\frac{(2k+1)\pi t}{2}\right)\Big],\,t\in(0,T).
\end{multline*}
Multiplying this last equality by $\cos\left(\frac{(2k+1)\pi t}{2}\right)$, for any $k\in\mathbb{N}$, and integrating on $(0,2)$ we obtain 
\begin{gather*}
 \label{balard}
a_k=0,\,k\in\mathbb{N},
\end{gather*}
and thus, 
\begin{gather}
\label{syrtes}
W^{\infty,1}\equiv 0\text{ on } (0,1).
\end{gather}

In particular, 
\begin{gather}
\label{syrtes2}
W^{\infty,1}(0)=\gamma_1\gamma_2 Z^{\infty}_3=0
\end{gather}
 and $$v^1_t(t,0)=\sum_{k=1}^{\infty}b_k\cos(\frac{(2k+1)\pi t}{2}),\,t\in(0,T).$$

 Parseval's identity then implies
 \begin{multline} \label{eq100}
 \int_0^T \vert v^1_t(t,0)\vert^2\,dt\ge \displaystyle\int_0^2\vert v^1_t(t,0)\vert^2\,dt= \sum\limits_{k=1}^{\infty} \vert b_k\vert^2
 =\vert W^{\infty,2}\vert^2_{L^2(0,1)}.
 \end{multline}
 Thus, from \eqref{eq80bis}, we finally get
$$W^{\infty,1}\equiv W^{\infty,2}\equiv0 \text{ on } (0,1)$$ 
and $$Z^{\infty}_1=Z^{\infty}_2=Z^{\infty}_3=0.$$
This is a contradiction with \eqref{th2eq3} and finishes the proof of Theorem \ref{th2}. \cqfd

%------------------------------------------------------------------------------------------------------------------------------
%------------------------------------------------------------------------------------------------------------------------------

\section{Numerical simulations}
\begin{figure}\label{fig1}
\begin{center}
\includegraphics[width=.8\textwidth]{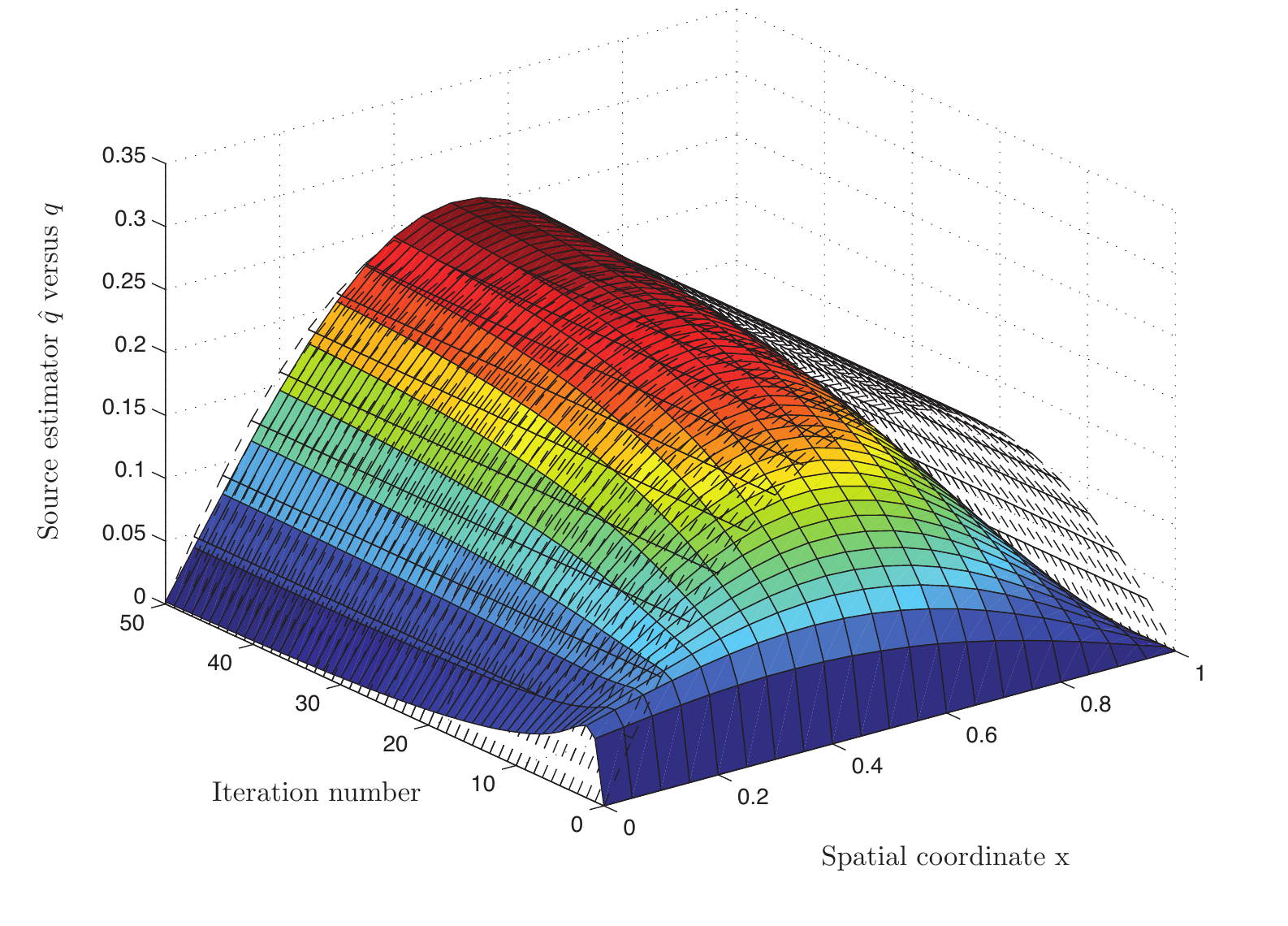}
\caption{The source estimator $\hat q$ is traced after each iteration; as it can be seen, the estimator after 50 iterations has converged towards $q$.}
\end{center}
\end{figure}
\begin{figure}\label{fig2}
\begin{center}
\includegraphics[width=.8\textwidth]{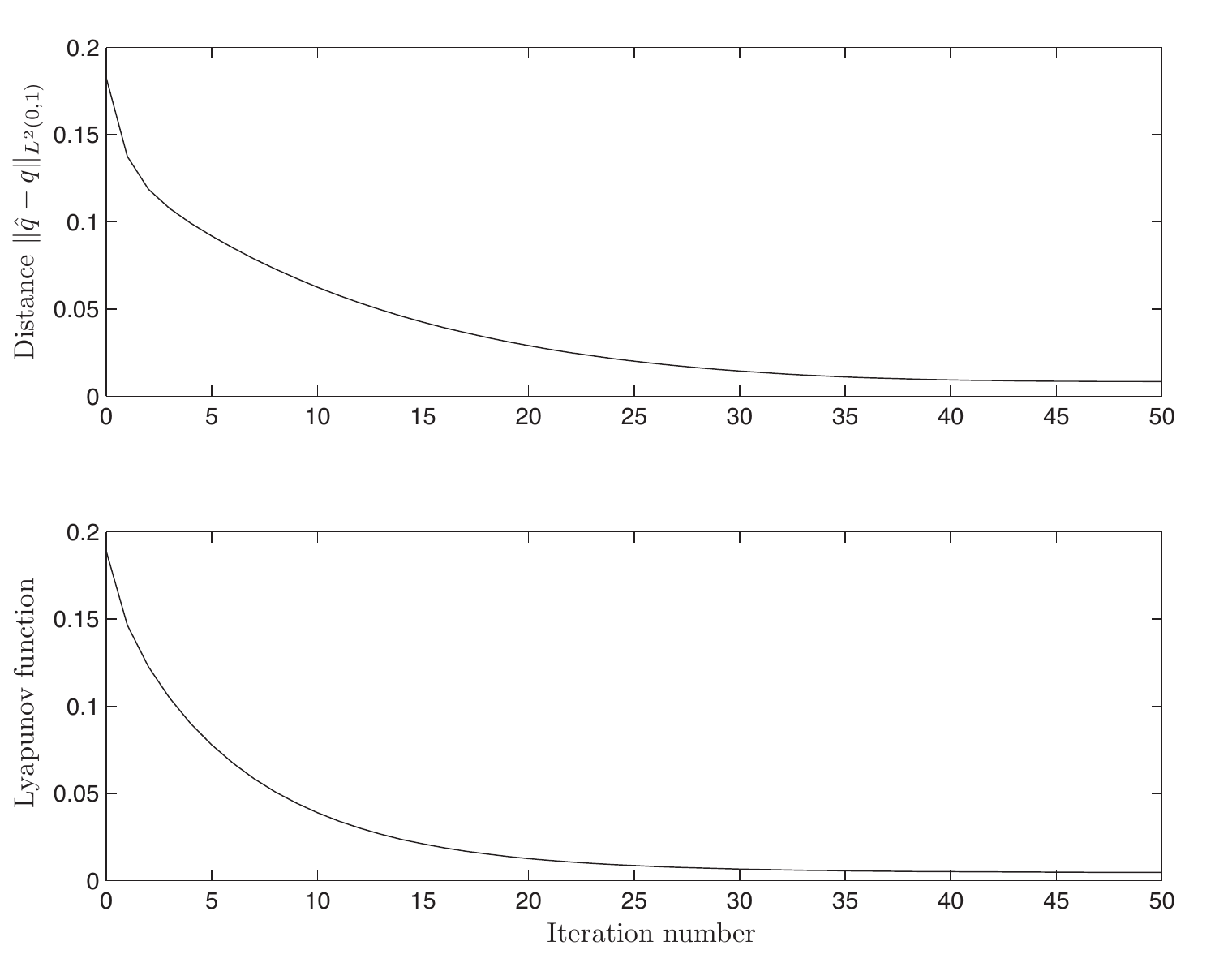}
\caption{The first plot illustrates the decrease of the $L^2$ distance between $\hat q$ and $q$ after each iteration and its convergence to zeros. The second plot illustrates the decrease of the Lyapunov function defined in~\eqref{eq82}.}
\end{center}
\end{figure}
In this section, we illustrate the efficiency of the above source estimation algorithm through numerical simulations. We consider the system~\eqref{eq1} with source term $q:=x-x^2$, together with the estimation algorithm~\eqref{eq13}-\eqref{eq14} with initial estimate $\hat{q}\equiv 0$. We fix the observation horizon to $T=3$ and consider 50 iterations of the estimator~\eqref{eq13}-\eqref{eq14}. The simulations of Figures~1and~2 illustrate the performance when we have added 10\% white noise on the measurement output and where the observer gains $\gamma_1$ and $\gamma_2$ are chosen to be
$$
\gamma_1=1,\qquad\gamma_2=1/2.
$$ 
The numerical simulations have been done through a finite difference method where the time and the space are discretized simultaneously. We have chosen a spatial discretization with 20 steps ($\Delta x=.05$) and a CFL coefficient of  $.005$ ($\Delta t=2.5e-04$). 

%------------------------------------------------------------------------------------------------------------------------------
%------------------------------------------------------------------------------------------------------------------------------

\section{The $N$-dimensional case}\label{sec:ext}

The aim of this section is to show that the Claim \ref{th1} can be easily extended to the $N$-dimensional case. More precisely, let $\Omega$ be a bounded domain in $\mathbb{R}^n$, with smooth boundary  $\Gamma$. We assume that $\Gamma$ is divided in two parts $\Gamma=\Gamma_0\cup \Gamma_1$, with $\Gamma_0\cap \Gamma_1=\emptyset$ and $\Gamma_0$ satisfying the geometrical optics condition~\cite{bardos-et-al-siam92}. That is, there exists a positive constant $L_{\Gamma_0}$  such that every generalized geodesic (the rays of geometrical optics) of length greater than $L_{\Gamma_0}$ passes through $\Gamma_0$ at a nondiffractive point.  Finally, let $T>0$,  $\omega\in\mathbb{R}$ and $q\in H^2(\Omega)\cap H^1_0(\Omega)$.  We denote $\mathcal{W}:=\{w\in H^1(\Omega) \text{ s.t. } w_{\vert \Gamma_1}=0\}.$ We consider the following system
\begin{equation}
\label{eq1n}
\left\{
\begin{array}{l}
v_{tt}-\Delta v=q(x)\cos(\omega t),\,(t,x)\in(0,T)\times \Omega,\\
v(t,x)=0,\,(t,x)\in(0,T)\times \Gamma,\\
v(0,x)=v_t(0,x)=0,\,x\in \Omega,\\
y(t,x)=\displaystyle\frac{\partial v}{\partial n}(t,x),\,(t,x)\in(0,T)\times \Gamma_0,
\end{array}
\right. 
\end{equation}
where $(u,u_t)$ represents the state of the system, $y$ is the output and we set $$\frac{\partial u}{\partial n}(t,x)=\sum\limits_{i=1}^{r}\nu_i(x)\frac{\partial u}{\partial x_i}(t,x),\,(t,x)\in(0,T)\times\Gamma,$$ where $\nu(x)=(\nu_1(x),...,\nu_r(x))$ is the outward unit normal to $\Gamma$ at $x$.\\
Then we have the following result
\begin{claim}
\label{th1n}
We can construct efficient observers for which the back and forth algorithm is convergent and which allow to reconstruct the unknown source term $q(x)$, using the measurement output $y(t)$ over the time-interval $(0,T)$, where $T>T_{obs}$ is any time strictly greater  than the minimal observability time of the wave equation observed at $\Gamma_0$~\cite{bardos-et-al-siam92}.
\end{claim}
To prove this claim, we proceed as in the proof of Theorem \ref{th1}.  We begin by introducing an equivalent system with $w=\omega^2 v+v_{tt}$, so that the unknown parameter is transferred to the initial state and we can easily see that the output of this second system, $Y$, is equivalent to the output of~\eqref{eq1n} in the space $H^2(0,T;L^2(\Gamma_0))$.
As before, we extend  $(w^1,w^2,z^1,z^2)$ to $((0,+\infty)\times\Omega)^2\times( (0,+\infty)\times\Gamma_0)^2$ and we propose exactly the same observer as in 1-dimensional case. The error dynamics are therefore given by

\small
\begin{equation}
\label{eq13n}
\left\{
\begin{array}{l}
\tilde{W}^1_{t}=\tilde{W}^2,\,(t,x)\in(2kT,(2k+1)T)\times\Omega,\vspace{0,1cm}\\
\tilde{W}^2_{t}=\Delta \tilde{W}^1,\,(t,x)\in(2kT,(2k+1)T)\times\Omega,\vspace{0,1cm}\\
\tilde{W}^1_{t}=-\tilde{W}^2,\,(t,x)\in((2k+1)T,(2k+2)T)\times\Omega,\vspace{0,1cm}\\
\tilde{W}^2_{t}=-\Delta \tilde{W}^1,\,(t,x)\in((2k+1)T,(2k+2)T)\times\Omega,\vspace{0,1cm}\\
\tilde{W}^1(t,x)=\gamma_1\tilde{Z}^1(t,x)+\gamma_1\gamma_2\tilde Z^3(t,x),\,(t,x)\in((2k+1)T,(2k+2)T)\times\Gamma_0,\vspace{0,1cm}
\end{array}
\right. 
\end{equation}
and
\begin{equation}
\label{eq14n}
\left\{
\begin{array}{l}
\tilde{Z}^1_t(t,x)=\tilde{Z}^2(t,x)-\gamma_2\tilde Z^1(t,x),\,(t,x)\in(2kT,(2k+1)T)\times\Gamma_0,\vspace{0,1cm}\\
\tilde{Z}^2_t(t,x)=\displaystyle-\omega^2\tilde{Z}^1(t,x)-\frac{\partial \tilde{W}^1}{\partial n}(t,x),\,(t,x)\in(2kT,(2k+1)T)\times\Gamma_0,\vspace{0,1cm}\\
\tilde{Z}^1_t(t,x)=-\tilde{Z}^2(t,x)-\gamma_2\tilde Z^1(t,x),\,(t,x)\in((2k+1)T,(2k+2)T)\times\Gamma_0,\vspace{0,1cm}\\
\tilde{Z}^2_t(t,x)=\displaystyle\omega^2\tilde{Z}^1(t,x)-\frac{\partial \tilde{W}^1}{\partial n}(t,x),\,(t,x)\in((2k+1)T,(2k+2)T)\times\Gamma_0,\vspace{0,1cm}\\
\tilde{Z}^3_t(t,x)=\tilde{Z}^1(t,x),\,(t,x)\in((2k+1)T,(2k+2)T)\times\Gamma_0,\vspace{0,1cm}
\end{array}
\right. 
\end{equation}
\normalsize

Here, the unique solution is defined on the state space
$$(\tilde W^1,\tilde W^2)\in C([0,+\infty);H^2(\Omega)\cap \mathcal{W})\times C([0,+\infty);\mathcal{W}),$$
$$(\tilde Z^1,\tilde Z^2,\tilde Z^3)\in H^1([0,+\infty);L^2(\Gamma_0))\times L^2([0,+\infty);L^2(\Gamma_0))\times H^2([0,+\infty);L^2(\Gamma_0))$$
Furthermore, we have the energy estimates
\small
\begin{multline}
\label{eq26n}
\vert \nabla \tilde W^1(t,.)\vert^2_{L^2(\Omega)}+\vert \tilde W^2(t,.)\vert^2_{L^2(\Omega)}+\gamma_1\vert \tilde Z^2(t,.)\vert^2_{L^2(\Gamma_0)}+\gamma_1\omega^2\vert \tilde Z^1(t,.)\vert^2_{L^2(\Gamma_0)}+2\gamma_1\gamma_2\omega^2\vert \tilde Z^1\vert^2_{L^2((0,t)\times \Gamma_0)}=\\
\vert \nabla \tilde W^1_0\vert^2_{L^2(\Omega)}+\vert \tilde W^2_0\vert^2_{L^2(\Omega)}+\gamma_1\vert \tilde Z^2_0\vert^2_{L^2(\Gamma_0)}+\gamma_1\omega^2\vert \tilde Z^1_0\vert^2_{L^2(\Gamma_0)}
\end{multline}\normalsize
and \small
\begin{multline}
\label{eq26bisn}
\frac{1}{2}(\vert \tilde W^2_{t}(t,.)\vert^2_{L^2(\Omega)}+\vert \nabla \tilde W^2(t,.)\vert^2_{L^2(\Omega)}+\gamma_1\omega^4\vert\tilde Z^1(t,.)\vert^2)_{L^2(\Gamma_0)}\\
+\frac{1}{4}\gamma_1\vert \frac{\partial \tilde W^1_x}{\partial n}(t,x)\vert^2_{L^2(\Gamma_0)} +\frac{1}{2}\gamma_1\gamma_2\omega^2\vert \tilde Z^2\vert^2_{L^2((0,t)\times \Gamma_0)}+\gamma_1\omega^2\vert \tilde Z^2(t,.)\vert^2_{L^2(\Gamma_0)} \\
 \le C_{\omega}(\vert \tilde W^1_0\vert^2_{H^2(\Omega)}+\vert \tilde W^2_0\vert^2_{H^1(\Omega)}+\vert \tilde Z^1_0\vert^2_{L^2(\Gamma_0)}+\vert \tilde Z^2_0\vert^2)_{L^2(\Gamma_0)},
\end{multline}\normalsize
where $C_{\omega}$ denotes a positive constant which only depends on $\omega$.

Finally,  assuming $T$ longer than the observability time of the wave equation, we have the following result,
\begin{gather}
\label{th2eq1n}
\lim\limits_{n\rightarrow +\infty}(\tilde{W}^1(2nT,.),\tilde{W}^2(2nT,.))= (0,0)\text{ in } \mathcal{W}\times L^2(\Omega),\\
\label{th2eq2n}
\lim\limits_{n\rightarrow +\infty}(\tilde Z^1(2nT,.),\tilde Z^2(2nT,.),\tilde Z3(2nT,.))= (0,0,0)\text{ in } L^2(\Gamma_0)^3.
\end{gather}

In order to prove~\eqref{th2eq1n} and~\eqref{th2eq2n}, one may proceed exactly as in the proof of Theorem \ref{th2}. The main difference is in the fact that by the energy inequalities, we have  $\tilde{Z}^1\in L^{\infty}(0,+\infty;L^2(\Gamma_0))$, and thus (up to a subsequence of $(\phi(n))_{n\ge0}$), there exists $Z^1_{\infty}\in L^2(\Gamma_0)$ such that 
\begin{gather}
\label{sentinelle}
\tilde{Z}^1(2\phi(n)T,.)\rightharpoonup Z^1_{\infty} \text{ weakly in } L^2(\Gamma_0)
\end{gather}
when $ {n\rightarrow +\infty}$, and $$\vert Z^1_{\infty}\vert_{L^2(\Gamma_0)}\le \liminf_{n\rightarrow +\infty} \vert \tilde{Z}^1(2\phi(n)T,.)\vert_{L^2(\Gamma_0)}.$$
We aim to prove that $\tilde{Z}^1(\phi(n)T,.)$ converges strongly to $0$ in $L^2(\Gamma_0)$
We start with  introducing the following Lyapunov function
\begin{eqnarray}
\begin{array}{l}
 \label{eq82n}
\mathcal{V}(\tilde{W}^1,\tilde{W}^2,\tilde{Z}^1,\tilde{Z}^2)(t):=\displaystyle\frac{1}{2}(\int_{\Omega}(\vert \nabla\tilde{W}^1(t,x)\vert^2+\vert \tilde{W}^2(t,x)\vert^2)\,dx\\
\phantom{nnnnnkkkkkkkkkkkkkkknnnn,nn}\displaystyle+\omega^2\int_{\Gamma_0}\vert \tilde{Z}^1(t,x)\vert^2\,dx+\int_{\Gamma_0}\vert \tilde{Z}^2(t,x)\vert^2\,dx),\,t\in(0,T).
\end{array}
\end{eqnarray}
Indeed, $\mathcal{V}(t)\ge0,\,t\ge0$ and applying \eqref{eq13n} and \eqref{eq14n}, one can compute, for any time $t\in(0,T)$,
\begin{gather}
 \label{eq82bisn}
\frac{d}{dt}\mathcal{V}(\tilde{W}^1,\tilde{W}^2,\tilde{Z}^1,\tilde{Z}^2)(t)=-\omega^2\int_{\Gamma_0}\vert \tilde{Z}^1(t,x)\vert^2\,dx\le 0.
\end{gather}
The function $\mathcal{V}$ is positive, decreasing. Consequently there exists $l\ge0$ such that
\begin{gather}
\label{eq88n}
\lim\limits_{t\rightarrow +\infty}\mathcal{V}(\tilde{W}^1,\tilde{W}^2,\tilde{Z}^1,\tilde{Z}^2)(t)=l.
\end{gather}

From \eqref{eq82bisn} and \eqref{eq88n}, 
\begin{gather}
\label{sentinelle2}
\displaystyle \int_0^{+\infty}\vert \tilde{Z}^1(t,.) \vert_{L^2(\Gamma_0)}^2\,dt=l-\mathcal{V}(0)<+\infty. 
\end{gather}
The idea is to apply the Barbalat's lemma. In this aim, we need to prove that $\frac{\partial}{\partial t} \vert \tilde{Z}^1(t,.) \vert_{L^2(\Gamma_0)}^2$ is bounded on $\mathbb{R}_+$ (which implies the absolute continuity).
This can be done through simple computations based on \eqref{eq13n} and \eqref{eq14n}. Consequently, from Barbalat's lemma, 
\begin{gather}
\label{sentinelle3}
\lim\limits_{t\rightarrow+\infty} \tilde{Z}^1(t,.)=0 \text{ in } L^2(\Gamma_0).
\end{gather}
In the same way, there exists $Z^2_{\infty}\in L^2(\Gamma_0)$ such that 
\begin{gather}
\label{sentinelle}
\tilde{Z}^2(2\phi(n)T,.)\rightharpoonup Z^2_{\infty} \text{ in } L^2(\Gamma_0)
\end{gather}
when $ {n\rightarrow +\infty}$. Once again, we want to prove that $Z^2_\infty$ is indeed 0 and furthermore that, the convergence is actually strong.  We first remark that $f(t):=\vert \tilde{Z}^2\vert^2_{L^2(\Gamma_0)}$ is absolutely continuous on $\mathbb{R}_+.$ Indeed, for any $t\in\mathbb{R}_+$,
\begin{eqnarray}
\begin{array}{l}
 \label{sentinelle4}
 f'(t)=-\displaystyle2\omega^2\int_{\Gamma_0}\tilde{Z}^2\tilde{Z}^1\,d\Gamma_0-2\int_{\Gamma_0}\tilde{Z}^2\frac{\partial \tilde{W}^1}{\partial n}\,d\Gamma_0,
\end{array}
\end{eqnarray}
and therefore, \eqref{eq26n} and \eqref{eq26bisn} imply the boundedness of $f'(t)$ on $\mathbb{R}_+$. 

Now, one has, at least formally, \small
\begin{eqnarray*}
\begin{array}{rcl}
 \label{sentinelle5}
\displaystyle\frac{d}{dt}(\int_{\Gamma_0}\tilde{Z}^1\tilde{Z}^2\,d\Gamma_0)&=&\displaystyle\int_{\Gamma_0}\tilde{Z}^1(-\omega^2\tilde{Z}^1-\frac{\partial \tilde{W}^1}{\partial n})\,d\Gamma_0+\int_{\Gamma_0}(\tilde{Z}^2-\gamma_2\tilde{Z}^1)\tilde{Z}^2\,d\Gamma_0\vspace{0,2cm}\\
&=&\displaystyle-\omega^2\int_{\Gamma_0}\vert\tilde{Z}^1\vert^2\,d\Gamma_0-\int_{\Gamma_0}\tilde{Z}^1\frac{\partial \tilde{W}^1}{\partial n}\,d\Gamma_0+\int_{\Gamma_0}\vert\tilde{Z}^2\vert^2\,d\Gamma_0-\gamma_2\int_{\Gamma_0}\tilde{Z}^1\tilde{Z}^2\,d\Gamma_0.
\end{array}
\end{eqnarray*}\normalsize
In other words,\small
\begin{eqnarray*}
\begin{array}{rcl}
 \label{sentinelle6}
\displaystyle -\int_{\Gamma_0}\tilde{Z}^1\frac{\partial \tilde{W}^1}{\partial n}\,d\Gamma_0+\int_{\Gamma_0}\vert\tilde{Z}^2\vert^2\,d\Gamma_0&=&\displaystyle\frac{d}{dt}(\int_{\Gamma_0}\tilde{Z}^1\tilde{Z}^2\,d\Gamma_0)+\omega^2\int_{\Gamma_0}\vert\tilde{Z}^1\vert^2\,d\Gamma_0+\gamma_2\int_{\Gamma_0}\tilde{Z}^1\tilde{Z}^2\,d\Gamma_0\vspace{0,1cm}\\
&\le&\displaystyle \displaystyle\frac{d}{dt}(\int_{\Gamma_0}\tilde{Z}^1\tilde{Z}^2\,d\Gamma_0)+(\omega^2+\frac{\gamma_2^2}{2})\int_{\Gamma_0}\vert\tilde{Z}^1\vert^2\,d\Gamma_0+\frac{1}{2}\int_{\Gamma_0}\vert\tilde{Z}^2\vert^2\,d\Gamma_0,
\end{array}
\end{eqnarray*}\normalsize
which finally gives
\begin{eqnarray*}
\begin{array}{rcl}
 \label{sentinelle7}
 \displaystyle g(t)+h(t)&\le& \displaystyle \displaystyle\frac{d}{dt}(\int_{\Gamma_0}\tilde{Z}^1\tilde{Z}^2\,d\Gamma_0)+(\omega^2+\frac{\gamma_2^2}{2})\int_{\Gamma_0}\vert\tilde{Z}^1\vert^2\,d\Gamma_0.
 \end{array}
\end{eqnarray*}
Here, $g(t):=-\int_{\Gamma_0}\tilde{Z}^1\frac{\partial \tilde{W}^1}{\partial n}\,d\Gamma_0$ which tends to $0$ when $t\rightarrow+\infty$ from  \eqref{eq26bisn} (which implies the boundedness in $L^2(\Gamma_0)$ of $\frac{\partial \tilde{W}^1}{\partial n}$ ) and \eqref{sentinelle3}. Also, $h(t):=\frac{1}{2}\int_{\Gamma_0}\vert\tilde{Z}^2\vert^2\,d\Gamma_0$, is positive and absolutely continuous from \eqref{sentinelle4}. Integrating the last inequality on $\mathbb{R}_+$ and applying \eqref{eq26n} (which implies the boundedness in $L^2(\Gamma_0)$ of $\tilde Z^2$),~\eqref{sentinelle3} and \eqref{sentinelle2} we obtain,
\begin{eqnarray*}
\begin{array}{rcl}
 \label{sentinelle8}
 \displaystyle \int_0^{+\infty}(g(t)+h(t))\,dt&\le&-\displaystyle \int_{\Gamma_0}\tilde{Z}^1(0,x)\tilde{Z}^2(0,x)\,d\Gamma_0+(\omega^2+\frac{\gamma^2_2}{2})\int_0^{+\infty}\int_{\Gamma_0}\vert\tilde{Z}^1\vert^2\,d\Gamma_0\,dt\\
 &<&+\infty.
\end{array}
\end{eqnarray*}
Barbalat's lemma implies that 
$$h(t)=\frac{1}{2}\int_{\Gamma_0}\vert\tilde{Z}^2\vert^2\,d\Gamma_0\rightarrow 0\qquad \text{as }t\rightarrow+\infty.$$ 
In particular, we have $$\tilde{Z}^2(2\phi(n)T,.)\rightarrow 0 \text{ in } L^2(\Gamma_0).$$ Now, from the boundedness of $\vert\tilde{Z}^1(2\phi(n)T,.)\vert_{L^2(\Gamma_0)}$ and $\vert\tilde{W}^1(2\phi(n)T,.)\vert_{L^2(\Gamma_0)}$ (which is a consequence of the convergence of $\tilde{W}^1(2\phi(n)T,.)$ in $\mathcal{W}$), and the fifth equation in~\eqref{eq13n}, we also have, up to a  subsequence, that there exists $Z^3_{\infty}\in H^{-1}(\Gamma_0)$ such that
\begin{gather}
\label{cvbisnn}
\lim\limits_{n\rightarrow +\infty}  \tilde{Z}^3(2\phi(n)T,.)=Z^3_{\infty} \text{ in } H^{-1}(\Gamma_0).
\end{gather}
As in the proof of Theorem \ref{th2}, we now define $(v^1,v^2)\in C([0,T];\mathcal{W})\times C([0,T];L^2(\Omega))$ and $(\zeta^1,\zeta^2,\zeta^3)\in H^1(0,T; L^2(\Gamma_0))\times L^2(0,T; L^2(\Gamma_0))\times H^2(0,T; H^{-1}(\Gamma_0))$ as the solution of
\begin{equation}
\label{eq80n}
\left\{
\begin{array}{l}
v^1_{t}=v^2,\,(t,x)\in(0,T)\times\Omega,\\
v^2_{t}=\Delta v^1,\,(t,x)\in(0,T)\times\Omega,\\
v^1=\gamma_1\zeta^1+\gamma_1\gamma_2\zeta^3,\,(t,x)\in(0,T)\times\Gamma_0,\\
\zeta^1_t=\zeta^2-\gamma_2 \zeta^1,\,(t,x)\in(0,T)\times\Gamma_0,\\
\zeta^2_t=-\omega^2\zeta^1(t)+\displaystyle \frac{\partial v^1}{\partial n},\,(t,x)\in(0,T)\times\Gamma_0,\\
\zeta^3_t=\zeta^1,\,(t,x)\in(0,T)\times\Gamma_0,\\
v^1(0,x)=W^1_{\infty}(x),\,v^2(0,x)=W^2_{\infty}(x),\,x\in\Omega,\\
\zeta^1(0,x)=0,\,\zeta^2(0,x)=0,\,\zeta^3(0)=Z^3_{\infty}(x),\,x\in\Gamma_0.
\end{array}
\right. 
\end{equation}
By continuity of flow with respect to the initial state, we have
\begin{gather}
 \label{eq84n}
\lim\limits_{n\rightarrow +\infty}(\tilde{W}^1(2\phi(n)T+t,.),\tilde{W}^2(2\phi(n)T+t,.))= (v^1,v^2)\text{ in } L^{\infty}((0,T);\mathcal{W})\times L^{\infty}((0,T);L^2(\Omega)),\\
\label{eq84bisn}
\lim\limits_{n\rightarrow +\infty}(\tilde{Z}^1(2\phi(n)T+t,.),\tilde{Z}^2(2\phi(n)T+t,.))=(\zeta^1,\zeta^2)\text{ in } L^{\infty}(0,T;L^2(\Gamma_0)),\\
\label{eq85bisn}
\lim\limits_{n\rightarrow +\infty}\tilde{Z}^3(2\phi(n)T+t,.)=\zeta^3\text{ in }L^{\infty}(0,T;H^{-1}(\Gamma_0).
\end{gather}

Following  the same strategy as in the proof of Theorem \ref{th2}, we have
$$
\zeta^1\equiv0, \quad \zeta^2\equiv0, \quad \zeta^3(t,x)= Z^3_{\infty}(x), \quad \frac{\partial v^1}{\partial n}\equiv0\qquad  \text{ on } (0,T)\times\Gamma_0
$$
Consequently, we obtain
\begin{equation}
\label{eq80bisn}
\left\{
\begin{array}{l}
v^1_{t}=v^2,\,(t,x)\in(0,T)\times\Omega,\\
v^2_{t}=\Delta v^1,\,(t,x)\in(0,T)\times\Omega,\\
v^1(t,x)=0,\,(t,x)\in(0,T)\times\Gamma_1,\\
\displaystyle \frac{\partial v^1}{\partial n}(t,x)=0,\,(t,x)\in(0,T)\times\Gamma_0,\\
v^1(t,x)=\gamma_1\gamma_2 Z^3_{\infty}(x),\,(t,x)\in(0,T)\times\Gamma_0,\\
v^1(0,x)=W^1_{\infty}(x),\,v^2(0,x)=W^2_{\infty}(x),\,x\in\Omega,
\end{array}
\right. 
\end{equation}
whose solutions are defined on the space $C([0,T];\mathcal{W})\times C([0,T];L^2(\Omega))$. As we have assumed $T>T_{\text{obs}}$, we can find a sufficiently small $\tau$ such that $T-\tau$ is also greater than $T_{\text{obs}}$: $T-\tau>T_{\text{obs}}$. We consider now the function
$$
(\psi^1(t,x),\psi^2(t,x)):=(v^1(t+\tau,x)-v^1(t,x),v^2(t+\tau,x)-v^2(t,x))
$$
defined on the space 
$C([0,T-\tau];\mathcal{W})\times C([0,T-\tau];L^2(\Omega))$. The functions $(\psi^1(t,x),\psi^2(t,x))$ satisfy the following equation:
\begin{equation}
\label{eq:final}
\left\{
\begin{array}{l}
\psi^1_{t}=\psi^2,\,(t,x)\in(0,T-\tau)\times\Omega,\\
\psi^2_{t}=\Delta \psi^1,\,(t,x)\in(0,T-\tau)\times\Omega,\\
\psi^1(t,x)=0,\,(t,x)\in(0,T-\tau)\times\Gamma,\\
\displaystyle \frac{\partial \psi^1}{\partial n}(t,x)=0,\,(t,x)\in(0,T)\times\Gamma_0,\\
\psi^1(0,x)=v^1(\tau)-W^1_{\infty}(x),\,\psi^2(0,x)=v^2(\tau)-W^2_{\infty}(x),\,x\in\Omega,
\end{array}
\right. 
\end{equation}
But as $T-\tau>T_{\text{obs}}$ this necessarily means that 
$$
\psi^1\equiv\psi^2\equiv 0. 
$$
As this is through for any $\tau$ such that $\tau\in [0,T_{\text{obs}}-T)$, we necessarily have that $v^1$ and $v^2$, solutions of~\eqref{eq80bisn}, are constants of time: 
$$
v^1(t,x)\equiv W^1_{\infty}(x),\qquad v^2(t,x)\equiv W^2_{\infty}(x).
$$ 
Inserting this inside~\eqref{eq80bisn}, $W^2_{\infty}(x)\equiv 0$ and $W^1_{\infty}(x)$ must satisfy the following elliptic equation: 
\begin{equation}
\label{eq:elliptic}
\left\{
\begin{array}{l}
\Delta W^1_{\infty} =0,\, x\in\Omega,\\
W^1_{\infty}(x)=0,\, x\in\Gamma_1,\\
W^1_{\infty}(x)=\gamma_1\gamma_2 Z^3_{\infty}(x),\,x\in\Gamma_0,\\
\displaystyle \frac{\partial W^1_{\infty}}{\partial n}(x)=0,\,x\in\Gamma_0,\\
\end{array}
\right. 
\end{equation}
The unique possible solution to this elliptic equation is trivially $W^1_\infty\equiv 0$ and $Z^3_\infty\equiv 0$. 
This ends the proof of the convergence and therefore implies the Claim~\ref{th1n}.\cqfd

\section{Appendix}

This appendix is devoted to the proof of Proposition \ref{prop2}. \\

\subsection{Existence part}
We first concentrate on the existence of a solution to the \textit{forward} system. More precisely, our aim is to prove the following result
\begin{proposition}
\label{prop7}
There exists $T_0>0$ such that for any $\gamma_1$, $\gamma_2>0$, for any $(q^0,q^1)\in H^2(0,1)\cap H^1_0(0,1)\times H^1_0(0,1)$ and  $(\xi_1^0,\xi_2^0,\xi_3^0)\in\mathbb{R}^3$, there exists $(v^1,v^2)\in C([0,T_0];L^2(0,1))\times C([0,T_0];H^{-1}(0,1))$ and $(\xi_1,\xi_2,\xi_3)\in H^1(0,T_0)\times L^2(0,T_0)\times H^2(0,T_0)$ solution to
\begin{equation}
\label{eq65}
\left\{
\begin{array}{l}
v^1_{t}=v^2,\,(t,x)\in(0,T_0)\times(0,1),\vspace{0,1cm}\\
v^2_{t}=v^1_{xx},\,(t,x)\in(0,T_0)\times(0,1),\vspace{0,1cm}\\
v^1(t,1)=0,\,v^1(t,0)=\gamma_1\xi_1(t)+\gamma_1\gamma_2\xi_3(t),\,t\in(0,T_0),\vspace{0,1cm}\\
v^1(0,x)=q^0(x),\,v^2(0,x)=q^1(x),\,x\in(0,1),
\end{array}
\right. 
\end{equation}
\begin{equation}
\label{eq65bis}
\left\{
\begin{array}{l}
\dot{\xi_1}(t)=\xi_2(t)-\gamma_2\xi_1(t),\,t\in(0,T_0),\vspace{0,1cm}\\
\dot{\xi_2}(t)=-\omega^2\xi_1(t)+v^1_x(t,0),\,t\in(0,T_0),\vspace{0,1cm}\\
\dot{\xi_3}(t)=\xi_1(t),\,t\in(0,T_0),\vspace{0,1cm}\\
\xi_1(0)=\xi_1^0,\, \xi_2(0)=\xi_2^0,\, \xi_2(0)=\xi_3^0.
\end{array}
\right. 
\end{equation}
\end{proposition}
\textbf{Proof of Proposition \ref{prop7}.}
Let $\gamma_1$, $\gamma_2>0$, $(q^0,q^1)\in H^2(0,1)\cap H^1_0(0,1)\times H^1_0(0,1)$ and  $(\xi_1^0,\xi_2^0,\xi_3^0)\in\mathbb{R}^3$. It follows from \cite[Th\'eor\`eme 4.2 page 46]{Lionsbook88} that for any $T>0$, for any $f\in H^1(0,T)$, there exists a unique solution $(v^1,v^2)$ to 
\begin{equation}
\label{eq27bis}
\left\{
\begin{array}{l}
v^1_{t}=v^2,\,(t,x)\in(0,T)\times(0,1),\\
v^2_{t}=v^1_{xx},\,(t,x)\in(0,T)\times(0,1),\\
v^1(t,1)=0,\,v^1(t,0)=f(t),\,t\in(0,T),\\
v^1(0,x)=q^0(x),\,v^2(0,x)=q^1(x),\,x\in(0,1),
\end{array}
\right. 
\end{equation}
with $(v^1,v^2)\in C([0,T];L^2(0,1))\times C([0,T];H^{-1}(0,1)).$ For any $T>0$ we define $\mathcal{F}_T$ on $H^1(0,T)$ by
\begin{gather}
\label{eq55}
 \mathcal{F}_T:f\in H^1(0,T)\mapsto \gamma_1\xi_1+\gamma_1\gamma_2\xi_3,
\end{gather}
where $\xi_1$ and $\xi_3$ satisfy
\begin{equation}
\label{eq56}
\left\{
\begin{array}{l}
\dot{\xi}_1(t)=\xi_2(t)-\gamma_2\xi_1(t),\,t\in(0,T),\\
\dot{\xi}_2(t)=-\omega^2\xi_1(t)+v^1_x(t,0),\,t\in(0,T),\\
\dot{\xi}_3(t)=\xi_1(t),\,t\in(0,T),\\
\xi_1(0)=\xi_1^0,\,\xi_2(0)=\xi_2^0,\,\xi_3(0)=\xi_3^0,
\end{array}
\right. 
\end{equation}
where $v^1\in C([0,T];L^2(0,1))$ denotes the solution of \eqref{eq27bis}. Our aim, here, is to prove that there exists some $T_0>0$ such that $\mathcal{F}_{T_0}$ has a fixed point. 

We first state and prove the following Lemma
\begin{lemma}
\label{lem3}
The solution $v^1$ to \eqref{eq27bis} moreover satisfies the following regularity property
\begin{gather}
\label{eq28}
\vert v^1_x(t,0)\vert^2_{L^2(0,T)}\le 2(4T^2+3)\vert f\vert^2_{H^1(0,T)}+2(2+T)(\vert q^0\vert^2_{H^1_0(0,1)}+\vert q^1\vert^2_{L^2(0,1)}).
\end{gather}
\end{lemma}
\textbf{Proof of Lemma \ref{lem3}.}
Let $t\in (0,T).$ We multiply the second equation in \eqref{eq27bis} by $v^2$ and integrate on $(0,t)\times(0,1).$ An integration by parts gives
\begin{eqnarray*}
\begin{array}{rcl}
\displaystyle\vert v^2(t,.)\vert^2_{L^2(0,1)}+\vert v^1_x(t,.)\vert^2_{L^2(0,1)}&=&\displaystyle-2\int_0^tv^2(s,0)v^1_x(s,0)\,ds+\vert q^1\vert^2_{L^2(0,1)}+\vert q_x^0\vert^2_{L^2(0,1)}\\
&\le& \displaystyle\epsilon_1\vert v^2(t,0)\vert^2_{L^2(0,T)}+\frac{1}{\epsilon_1}\vert v^1_x(t,0)\vert^2_{L^2(0,T)}\vspace{0,1cm}\\
&&\displaystyle\phantom{ttttnnnnnnnnttttt}+\vert q^1\vert^2_{L^2(0,1)}+\vert q_x^0\vert^2_{L^2(0,1)},
\end{array}
\end{eqnarray*}
where $\epsilon_1$ denotes a positive constant that will be chosen later.
This last inequlity holding for any time $t\in(0,T)$, we finally obtain
\begin{eqnarray}
\begin{array}{l}
\label{eq50}
\displaystyle \sup_{t\in (0,T)}(\vert v^2(t,.)\vert^2_{L^2(0,1)}+\vert v^1_x(t,.)\vert^2_{L^2(0,1)})\le \displaystyle\epsilon_1\vert v^2(t,0)\vert^2_{L^2(0,T)}\vspace{0,1cm}\\
\displaystyle\phantom{nlljjjmllllllmmmmmmmkkkkjjjjjt}+\frac{1}{\epsilon_1}\vert v^1_x(t,0)\vert^2_{L^2(0,T)}+\vert q^1\vert^2_{L^2(0,1)}+\vert q^0_x\vert^2_{L^2(0,1)}.
\end{array}
\end{eqnarray}
Now, we multiply the second equation in \eqref{eq27bis} by $(x-1)v^1_x$ and integrate on $(0,T)\times(0,1)$. Several integrations by parts  give
\begin{eqnarray*}
\begin{array}{l}
\displaystyle
\frac{1}{2}(\int_0^T\int_0^1\vert v^2(t,x)\vert^2\,dx\,dt-\int_0^T\vert v^2(t,0)\vert^2\,dt+\int_0^T\int_0^1\vert v^1_x(t,x)\vert^2\,dx\,dt)\vspace{0,1cm}\\
\displaystyle\phantom{tttt}+\int_0^1(x-1)v^1_x(T,x)v^1_t(T,x)\,dx-\int_0^1(x-1)v^1_x(0,x)v^1_t(0,x)\,dx-\frac{1}{2}\int_0^T\vert v^1_x(t,0)\vert^2\,dt=0.
\end{array}
\end{eqnarray*}
Using the boundary conditions in \eqref{eq27bis} this leads
\begin{eqnarray}
\begin{array}{l}
\label{eq51}
\displaystyle
\vert v^1_x(t,0)\vert^2_{L^2(0,T)}\le\displaystyle \int_0^T\int_0^1(\vert v^2(t,x)\vert^2\,dx\,dt+\vert v^1_x(t,x)\vert^2)\,dx\,dt-\vert \dot f\vert^2_{L^2(0,T)}\vspace{0,1cm}\\
\displaystyle \phantom{tttttttttttttttttttttttttt}+\vert v^1_x(T,.)\vert^2_{L^2(0,1)}+\vert v^2(T,.)\vert^2_{L^2(0,1)}+\vert q^0_x\vert^2_{L^2(0,1)}+\vert q^1\vert^2_{L^2(0,1)}.
\end{array}
\end{eqnarray}
Combining now \eqref{eq50} and \eqref{eq51} and taking $\epsilon_1:= 4T$, we obtain
\begin{eqnarray*}
\begin{array}{l}
\label{eq52}
\displaystyle
\frac{3}{4}\vert v^1_x(t,0)\vert^2_{L^2(0,T)}\le (4T^2-1)\vert f\vert^2_{H^1(0,T)}+\vert v^1_x(T,.)\vert^2_{L^2(0,1)}+\vert v^2(T,.)\vert^2_{L^2(0,1)}\\
\displaystyle \phantom{tttttttttlllllmmmmmmmmmlllllllllttttttttttttttttt}+(1+T)(\vert q^0_x\vert^2_{L^2(0,1)}+\vert q^1\vert^2_{L^2(0,1)}).
\end{array}
\end{eqnarray*}
We use one more time \eqref{eq50} with $\epsilon_1:=4$ and we finally get
\begin{eqnarray*}
\begin{array}{l}
\label{eq53}
\displaystyle
\vert v^1_x(t,0)\vert^2_{L^2(0,T)}\le 2(4T^2+3)\vert f\vert^2_{H^1(0,T)}+2(2+T)(\vert q^0_x\vert^2_{L^2(0,1)}+\vert q^1\vert^2_{L^2(0,1)}),
\end{array}
\end{eqnarray*}
which ends the proof of Lemma \ref{lem3}. 
\cqfd

Now, we prove the following proposition ensuring that the map $\mathcal{F}_T$ is well-defined. 
\begin{proposition}
\label{prop5}
For any $T>0$, $\mathcal{F}_T:H^1(0,T)\rightarrow H^1(0,T)$ is well-defined. More precisely, for any $f\in H^1(0,T)$ there exists a unique solution $(\xi_1,\xi_2,\xi_3)\in H^1(0,T)\times L^2(0,T)\times H^2(0,T)$ to \eqref{eq56}.
\end{proposition}
\textbf{Proof of Proposition \ref{prop5}.} Let $T>0$ and $f\in H^1(0,T)$. We define $A\in \mathcal{M}_2(\mathbb{R})$ and $B\in \mathcal{M}_{2,1}(\mathbb{R})$ as
\begin{gather*}
A:=\left(\begin{array}{cc}-\gamma_2 & 1 \\-\omega^2 & 0\end{array}\right),\qquad
B:=\left(\begin{array}{c}0 \\v^1_x(t,0)\end{array}\right),
\end{gather*}
where $(v^1,v^2)$ denotes the solution of \eqref{eq27bis} associated to $f$.
It is clear that $(\xi_1,\xi_2,\xi_3)$ defined by,
\begin{gather}
\label{bisbald}
\left(\begin{array}{c}\xi_1(t) \\\xi_2(t)\end{array}\right):=e^{tA}\left(\begin{array}{c}\xi_1^0 \\\xi_2^0\end{array}\right)+\int_0^te^{(t-s)A}B(s)\,ds
\quad \text{ and }\quad \xi_3(t):=\xi_3^0+\int_0^t\xi_1(s)\,ds
\end{gather}
is the unique solution of \eqref{eq56}. Let us prove that $(\xi_1,\xi_2,\xi_3)\in H^1(0,T)\times L^2(0,T)\times H^2(0,T)$. Simple computations lead to 
\begin{eqnarray}
\begin{array}{rcl}
\label{eq57}
\vert (\xi_1,\xi_2)\vert^2_{L^2(0,T)}&\le& \displaystyle\vert e^{tA}\left(\begin{array}{c}\xi_1^0 \\\xi_2^0\end{array}\right)\vert^2_{L^2(0,T)}+\vert \int_0^te^{(t-s)A}B(s)\,ds\vert_{L^2(0,T)}\vspace{0,1cm}\\
&\le&  \displaystyle C T e^{TC}(\vert\xi_1^0\vert^2+\vert\xi_2^0\vert^2)+C\int_0^T  t e^{Ct}\left(\int_0^t |v^1_x(s,0)|^2\,ds\right)\,dt,
\end{array}
\end{eqnarray}
where we recall that $C=C(\omega,\gamma_1,\gamma_2)$ denotes various positive constants which may vary from line to line and only depend on $\omega$, $\gamma_1$ and $\gamma_2$.
%%%%%%%%%%%%%%%%%%%%%%%%%%%
Applying Lemma \ref{lem3}, \eqref{eq57} becomes
\begin{eqnarray}
\begin{array}{rcl}
\label{eq59}
\vert (\xi_1,\xi_2)\vert^2_{L^2(0,T)}&\le& \displaystyle CT e^{TC}(\vert\xi_1^0\vert^2+\vert\xi_2^0\vert^2)+CT^2e^{CT}\vert v^1_x(t,0)\vert^2_{L^2(0,T)}\vspace{0,1cm}\\
&\le& \displaystyle CT e^{TC}(\vert\xi_1^0\vert^2+\vert\xi_2^0\vert^2)\vspace{0,1cm}\\
&& \displaystyle\phantom{tttttt}+2Ce^{CT}T^2\left((4T^2+3)\vert f\vert^2_{H^1(0,T)}+(2+T)(\vert q^0\vert^2_{H^1_0(0,1)}+\vert q^1\vert^2_{L^2(0,1)})\right).
\end{array}
\end{eqnarray}
From the last inequality, $(\xi_1,\xi_2)\in L^2(0,T)$. Consequently, from the first line in \eqref{eq56}, 
\begin{gather}
\label{eq58}
\xi_1\in H^1(0,T).
\end{gather}
 The fact that $\xi_3\in H^2(0,T)$ finally comes from \eqref{bisbald} and \eqref{eq58}. This ends the proof of Proposition \ref{prop5}. \cqfd
 %%%%%%%%%%%%
 Now, we prove the following result, implying the existence part of the Proposition \ref{prop2}
 \begin{proposition}
 \label{prop6}
 There exists sufficiently small $T_0$ only depending on $\omega$, $\gamma_1$ and $\gamma_2$, such that $\mathcal{F}_{T_0}$ admits a fixed point in $H^1(0,T_0)$. 
 \end{proposition}
\textbf{Proof of Proposition \ref{prop6}.} 
For any $R>0$ we call $B_R:=\{f\in H^1(0,T) \text{ s. t. } \vert f\vert_{H^1(0,T)}\le R\}$.  In order to prove Proposition \ref{prop6}, we are going to prove the two following points, for sufficiently small $T_0$, \\
i) there exists $R_0>0$ such that $\mathcal{F}_{T_0}(B_{R_0})\subset B_{R_0}$,\\
ii) $\mathcal{F}_{T_0}$  is contractant.\\
\underline{Point i):}\vspace{0,1cm}\\
Let $R_0>0$ be such that 
\begin{gather}
\label{eq61}
R_0>\vert\xi_1^0\vert^2+\vert\xi_2^0\vert^2+\vert\xi_3^0\vert^2+\vert q^0\vert_{H^2(0,1)}^2+\vert q^1\vert^2_{H^1_0(0,1)},
\end{gather}
and let $f\in B_{R_0}.$ From \eqref{eq56}, \eqref{bisbald} and \eqref{eq59}, 
\begin{eqnarray}
\begin{array}{rcl}
\label{eq60}
\vert \xi_3\vert_{H^1(0,T)}^2&\le&\vert \xi_3\vert^2_{L^2(0,T)}+ \vert\xi_1\vert^2_{L^2(0,T)}\vspace{0,1cm}\\
&\le& T\vert \xi_3^0\vert^2+T^2\vert \xi_1\vert^2_{L^2(0,T)}+ \vert \xi_1\vert^2_{L^2(0,T)}\\
&\le& T\vert \xi_3^0\vert^2+C(1+T^2)( \displaystyle T e^{TC}(\vert\xi_1^0\vert^2+\vert\xi_2^0\vert^2)\vspace{0,1cm}\\
&& \displaystyle\phantom{tttttttt}+2Ce^{CT}T^2\left((4T^2+3)\vert f\vert^2_{H^1(0,T)}+(2+T)(\vert q^0_x\vert^2_{L^2(0,1)}+\vert q^1\vert^2_{L^2(0,1)}))\right).
\end{array}
\end{eqnarray}
By straightforward computations, there exists sufficiently small $T_0$, only depending on $(\omega,\gamma_1,\gamma_2)$, such that 
\begin{gather*}
\vert \gamma_1 \xi_1+\gamma_1\gamma_2\xi_3\vert_{H^1(0,T_0)}<R_0,
\end{gather*}
which implies the Point i).\\
\underline{Point ii):}\\
Let $f_1$, $f_2\in B_{R_0}.$ We denote by $(v^1,v^2)$ and $(\xi_1,\xi_2,\xi_3)$ (respectively by $(\eta^1,\eta^2)$ and $(\zeta_1,\zeta_2,\zeta_3)$) the solutions to \eqref{eq27bis} and \eqref{eq56} associated with $f_1$ (respectively with $f_2$). Then  $(\tilde{v}^1,\tilde{v}^2):=(v^1-\eta^1,v^2-\eta^2)$ and $(\tilde{\xi}_1,\tilde{\xi}_2,\tilde{\xi}_3):=(\xi_1,\xi_2,\xi_3)-(\zeta_1,\zeta_2,\zeta_3)$ satisfy
\begin{equation}
\label{eq62}
\left\{
\begin{array}{l}
\tilde{v}_{t}^1=\tilde{v}^2,\,(t,x)\in(0,T)\times(0,1),\\
\tilde{v}_{t}^2=\tilde{v}^1_{xx},\,(t,x)\in(0,T)\times(0,1),\\
\tilde{v}^1(t,1)=0,\,\tilde{v}^1(t,0)=f_1(t)-f_2(t),\,t\in(0,T),\\
\tilde{v}^1(0,x)=0,\,\tilde{v}^2(0,x)=0,\,x\in(0,1),
\end{array}
\right. 
\end{equation}
and
\begin{equation}
\label{eq63}
\left\{
\begin{array}{l}
\dot{\tilde{\xi}}_1(t)=\tilde{\xi}_2(t)-\gamma_2\tilde{\xi}_1(t),\,t\in(0,T),\\
\dot{\tilde{\xi}}_2(t)=-\omega^2\tilde{\xi}_1(t)+\tilde{v}^1_x(t,0),\,t\in(0,T),\\
\dot{\tilde{\xi}}_3(t)=\tilde{\xi}_1(t),\,t\in(0,T),\\
\tilde{\xi}_1(0)=0,\,\tilde{\xi}_2(0)=0,\,\tilde{\xi}_3(0)=0.
\end{array}
\right. 
\end{equation}
Adapting Lemma \ref{lem3}, \eqref{eq59} and \eqref{eq60}, one easily obtains 
\begin{eqnarray*}
\begin{array}{rcl}
\label{eq64}
\vert \gamma_1\tilde{\xi}_1+\gamma_1\gamma_2\tilde{\xi}_3\vert^2_{H^1(0,T)}&\le&\widetilde{C}e^{CT}T^2(4T^2+3)\vert f_1-f_2\vert^2_{H^1(0,T)},
\end{array}
\end{eqnarray*}
where $\widetilde{C}$ is a positive constant only depending on $(\omega,\gamma_1,\gamma_2)$. Therefore for sufficiently small $T_0$ only depending on $(\omega,\gamma_1,\gamma_2)$, $\widetilde{C}e^{CT_0}T_0^2(4T_0^2+3)<1$ and the Point ii is proved.\\
From Point i) and Point ii), $\mathcal{F}_{T_0}$ admits a fixed point in $H^1(0,T_0)$: this ends the proof of Proposition \ref{prop6}. We emphasize the fact that the time $T_0$ defined in this proposition, does not depend on the initial conditions $(q^0,q^1)\in H^2(0,1)\cap H^1_0(0,1)\times H^1_0(0,1)$ and $(\xi_1^0,\xi_2^0,\xi_3^0)\in\mathbb{R}^3$. \cqfd
In the same way, one obtains the following result (existence of a solution to the system \textit{backward}
\begin{proposition}
\label{prop7bis}
There exists $T_1>0$ such that for any $\gamma_1$, $\gamma_2>0$, for any $(q^0,q^1)\in H^2(0,1)\cap H^1_0(0,1)\times H^1_0(0,1)$ and  $(\xi_1^0,\xi_2^0,\xi_3^0)\in\mathbb{R}^3$, there exists $(v^1,v^2)\in C([0,T_1];L^2(0,1))\times C([0,T_1];H^{-1}(0,1))$ and $(\xi_1,\xi_2,\xi_3)\in H^1(0,T_1)\times L^2(0,T_1)\times H^2(0,T_1)$ solution to
\begin{equation}
\label{eq65*}
\left\{
\begin{array}{l}
v^1_{t}=-v^2,\,(t,x)\in(0,T_1)\times(0,1),\vspace{0,1cm}\\
v^2_{t}=-v^1_{xx},\,(t,x)\in(0,T_1)\times(0,1),\vspace{0,1cm}\\
v^1(t,1)=0,\,v^1(t,0)=\gamma_1\xi_1(t)+\gamma_1\gamma_2\xi_3(t),\,t\in(0,T_1),\vspace{0,1cm}\\
v^1(0,x)=q^0(x),\,v^2(0,x)=q^1(x),\,x\in(0,1),
\end{array}
\right. 
\end{equation}
\begin{equation}
\label{eq65bis*}
\left\{
\begin{array}{l}
\dot{\xi_1}(t)=-\xi_2(t)-\gamma_2\xi_1(t),\,t\in(0,T_1),\vspace{0,1cm}\\
\dot{\xi_2}(t)=\omega^2\xi_1(t)-v^1_x(t,0),\,t\in(0,T_1),\vspace{0,1cm}\\
\dot{\xi_3}(t)=\xi_1(t),\,t\in(0,T_1),\vspace{0,1cm}\\
\xi_1(0)=\xi_1^0,\, \xi_2(0)=\xi_2^0,\, \xi_2(0)=\xi_3^0.
\end{array}
\right. 
\end{equation}
\end{proposition}
We are now going to prove that  times $T_0$ and $T_1$ defined in Proposition \ref{prop7} and \ref{prop7bis} can actually be chosen arbitrarily.

%%%%%%%%%%%%%%%%%%%%%%%%%%%%%%%%%%%%%%%%
%%%%%%%%%%%%%%%%%%%%%%%%%%%%%%%%%%%%%%%

\subsection{Global well-posedness}

%%%%%%%%%%%%%%%%%%%%%%%%%%%%%%%%%%%%%%%%%%%
In the context of an established local well-posedness result, it suffices to prove the following global a priori estimate for smooth solutions to system  \eqref{eq65}-\eqref{eq65bis} (the fact that in this case it will also hold for the  solution of system \eqref{eq65*}-\eqref{eq65bis*} is obvious)
\begin{eqnarray}
\begin{array}{l}
\label{eq26*}
\vert v^1_x(T_0,.)\vert^2_{L^2(0,1)}+\vert v^2(T_0,.)\vert^2_{L^2(0,1)}+\gamma_1\vert\xi_2(T_0)\vert^2+\gamma_1\omega^2\vert\xi_1(T_0)\vert^2+2\gamma_1\gamma_2\omega^2\vert\xi_1\vert^2_{L^2(0,T_0)}=\vspace{0,2cm}\\
\phantom{ttttttttttttttttlllllllllllllllllllbbbbbbbbbbllllllllttttttttt}\vert q^0_x\vert^2_{L^2(0,1)}+\vert q^1\vert^2_{L^2(0,1)}+\gamma_1\vert\xi_2^0\vert^2+\gamma_1\omega^2\vert\xi_1^0\vert^2,
\end{array}
\end{eqnarray}
and 
\begin{eqnarray}
\begin{array}{l}
\label{eq26bis*}
\displaystyle \frac{1}{2}(\vert v^2_{t}(T_0,.)\vert^2_{L^2(0,1)}+\vert v^2_{x}(T_0,.)\vert^2_{L^2(0,1)}+\gamma_1\omega^4\vert\xi_1(T_0)\vert^2)+\frac{1}{4}\gamma_1\vert v^1_x(T_0,0)\vert^2\\
\displaystyle\phantom{gglllllkkkkg}
+\frac{1}{2}\gamma_1\gamma_2\omega^2\int_0^{T_0}\vert \xi_2(s)\vert^2\,ds+\gamma_1\omega^2\vert \xi_2(T_0)\vert^2\le \\
\displaystyle \phantom{ggkkkkkkkkkllllllllllllllllllllllllllllkkg} C_{\omega,\gamma_1,\gamma_2}(\vert q^0\vert^2_{H^2(0,1)}+\vert q^1\vert^2_{H^1(0,1)}+\vert \xi_1^0\vert^2+\vert \xi_2^0\vert^2),
\end{array}
\end{eqnarray}
where $C_{\omega,\gamma_1,\gamma_2}$ denotes a positive constant which only depends on $\omega$, $\gamma_1$ and $\gamma_2$. One can note that these inequalities also ensure more regularity for solutions to \eqref{eq65}-\eqref{eq65bis}.\\

\underline{Step 1: Proof of \eqref{eq26*}.}\vspace{0,1cm}

Let $T_0$ be defined as in the Proposition~\ref{prop7} and let $(v^1,v^2,\xi_1,\xi_2,\xi_3)\in C([0,T_0];L^2(0,1))\times C([0,T_0];H^{-1}(0,1))\times H^1(0,T_0)\times L^2(0,T_0)\times H^2(0,T_0)$ denote a solution to \eqref{eq65}. We first multiply the second equation in \eqref{eq65} by $v^2$ and integrate on $(0,1)$. A straightforward integration by parts gives 
\begin{gather*}
\frac{1}{2}\partial_t(\int_0^1\vert v^2\vert^2+\vert v^1_x\vert^2)+v^1_x(t,0)v^2(t,0)=0.
\end{gather*}

This equality implies, from the boundary conditions in \eqref{eq65}
\begin{gather*}
\frac{1}{2}\partial_t(\int_0^1\vert v^2\vert^2+\vert v^1_x\vert^2)+(\dot{\xi_2}+\omega^2\xi_1)(\gamma_1\dot{\xi_1}+\gamma_1\gamma_2\dot{\xi_3})=0,
\end{gather*}
%\begin{eqnarray}
%\begin{array}{}
hence
\begin{gather*}
\frac{1}{2}\partial_t(\int_0^1(\vert v^2\vert^2+\vert v^1_x\vert^2)\,dx+\gamma_1\vert\xi_2\vert^2+\gamma_1\omega^2\vert \xi_1\vert^2)+\gamma_1\gamma_2\omega^2\vert\xi_1\vert^2=0
\end{gather*}
and \eqref{eq26*} follows. (We first obtain the previous equality in the case where $(v^1,v^2)$ is assumed to be smooth enough and we then use a density argument to obtain it in the general case).\\

\underline{Step 2: Proof of \eqref{eq26bis*}.}\vspace{0,1cm}

Now we derive the second equation in \eqref{eq65} with respect to the time variable and we integrate on $(0,1)$. (In the same manner, we first assume that $(v^1,v^2)$ is smooth enough and we then use a density argument to obtain it in the general case.)  It follows,
\begin{gather}
\label{eq41}
\frac{1}{2}\partial_t\int_0^1(\vert v^2_{t}\vert^2+\vert v^2_{x}\vert^2)+v^2_{x}(t,0)v^2_{t}(t,0)=0.
\end{gather}
Using the boundary conditions in \eqref{eq65} we compute,
\begin{eqnarray*}
\begin{array}{rcl}
\label{eq40}
v^2_{x}(t,0)v^2_{t}(t,0)&=&v^2_{x}(t,0)\gamma_1\dot{\xi_2}\\
&=&\gamma_1v^2_{x}(t,0)(-\omega^2\xi_1+v^1_x(t,0))\\
&=&-\gamma_1\omega^2 v^2_{x}(t,0)\xi_1+\frac{1}{2}\gamma_1\partial_t(\vert v^1_x(t,0)\vert^2).
\end{array}
\end{eqnarray*}
We integrate \eqref{eq41} on $(0,T_0)$ and use the initial conditions in \eqref{eq65} and the last identity to obtain
\begin{eqnarray}
\begin{array}{l}
\label{eq42}
\displaystyle \frac{1}{2}(\vert v^2_{t}(T_0,.)\vert^2_{L^2(0,1)}+\vert v^2_{x}(T_0,.)\vert^2_{L^2(0,1)}+\gamma_1\vert v^1_x(T_0,0)\vert^2)-\gamma_1\omega^2\int_0^{T_0}v^2_{x}(t,0)\xi_1(t)\,dt=\\
\displaystyle\phantom{ttttttlllllllmmmmmmmmmmmmmmmlltt}\frac{1}{2}(\vert q^0_{xx}\vert^2_{L^2(0,1)}+\vert q^1_x\vert^2_{L^2(0,1)}+\gamma_1\vert q^0_x(0)\vert^2).
\end{array}
\end{eqnarray}
Using one more time the boundary conditions in \eqref{eq65}, several integrations by parts give
\begin{eqnarray}
\begin{array}{rcl}
\label{eq43}
\displaystyle-\gamma_1\omega^2\int_0^{T_0}\xi_1(t)v^2_{x}(t,0)\,dt&=&\displaystyle\gamma_1\omega^2\int_0^{T_0}\dot{\xi_1}(t)v^1_x(t,0)\,dt-\gamma_1\omega^2[\xi_1(t)v^1_x(t,0)]^{T_0}_0\\
&=&\displaystyle\gamma_1\omega^2\int_0^{T_0}\dot{\xi_1}(t)(\dot{\xi_2}(t)+\omega^2\xi_1(t))dt-\gamma_1\omega^2[\xi_1(t)v^1_x(t,0)]^{T_0}_0\vspace{0,1cm}\\
&=&\displaystyle\frac{\gamma_1\omega^4}{2}(\vert\xi_1({T_0})\vert^2-\vert\xi_1^0\vert^2)-\gamma_1\omega^2\int_0^{T_0}\ddot{\xi_1}(t)\xi_2(t)\,dt\vspace{0,1cm}\\
&&\displaystyle\phantom{tkkkkktttt}+\gamma_1\omega^2[\dot{\xi_1}(t)\xi_2(t)]_0^{T_0}-\gamma_1\omega^2[\xi_1(t)v^1_x(t,0)]^{T_0}_0.
\end{array}
\end{eqnarray}
Now,
\begin{eqnarray}
\begin{array}{rcl}
\label{eq44}
\displaystyle -\gamma_1\omega^2\int_0^{T_0}\ddot{\xi_1}(t)\xi_2(t)\,dt&=&\displaystyle-\gamma_1\omega^2\int_0^{T_0}(\dot{\xi_2}(t)-\gamma_2\xi_2(t)+\gamma_2^2\xi_1(t))\xi_2(t)\,dt\\
&=&\displaystyle-\frac{\gamma_1\omega^2}{2}(\vert\xi_2({T_0})\vert^2-\vert \xi_2^0\vert^2)+\gamma_1\gamma_2\omega^2\int_0^{T_0}\vert\xi_2(t)\vert^2\,dt\vspace{0,1cm}\\
&&\phantom{mmmmmmmmmmmmmmm}\displaystyle-\gamma_1\gamma_2^2\omega^2\int_0^{T_0}\xi_1(t)\xi_2(t)\,dt.
\end{array}
\end{eqnarray}
Combining now \eqref{eq42}-\eqref{eq44}, we obtain
\begin{eqnarray}
\begin{array}{l}
\label{eq45}
\displaystyle \frac{1}{2}(\vert v^2_{t}(T_0,.)\vert^2_{L^2(0,1)}+\vert v^2_{x}(T_0,.)\vert^2_{L^2(0,1)}+\gamma_1\vert v^1_x(T_0,0)\vert^2+\gamma_1\omega^4\vert\xi_1(T_0)\vert^2+\gamma_1\omega^2\vert\xi_2^0\vert^2)\\
\displaystyle\phantom{ggkkkkg}+\gamma_1\gamma_2\omega^2\int_0^{T_0}\vert \xi_2(t)\vert^2\,dt+\gamma_1\omega^2\vert \xi_2(T_0)\vert^2-\gamma_1\gamma_2\omega^2\xi_1(T_0)\xi_2(T_0)+\gamma_1\omega^2\xi_1^0q^0_x(0)=\\
\displaystyle \phantom{ggkoooooooooookkkkkkg} \frac{1}{2}(\vert q^0_{xx}\vert^2_{L^2(0,1)}+\vert q^1_x\vert^2_{L^2(0,1)}+\gamma_1\vert q^0_x(0)\vert^2+\gamma_1\omega^4\vert\xi_1^0\vert^2+\gamma_1\omega^2\vert \xi_2(T_0)\vert^2)\vspace{0,1cm}\\
\displaystyle\phantom{ggkkkkpmmpkkkg}+\gamma_1\gamma_2^2\omega^2\int_0^{T_0}\xi_1(t)\xi_2(t)\,dt+\omega^2\vert\xi_2^0\vert^2-\gamma_1\gamma_2\omega^2\xi_1^0\xi_2^0+\gamma_1\omega^2\xi_1(T_0)v^1_x(T_0,0).
\end{array}
\end{eqnarray}
Using the following estimates
\begin{gather*}
\gamma_1\gamma_2\omega^2\xi_1(T_0)\xi_2(T_0)\le \displaystyle \frac{\gamma_1\gamma_2\omega^2}{2}(\vert \xi_1(T_0)\vert^2+\vert \xi_2(T_0)\vert^2),\\
\gamma_1\omega^2\xi_1^0q_x^0(0)\le \displaystyle \frac{\gamma_1\omega^2}{2}(\vert \xi_1^0\vert^2+\vert q_x^0(0)\vert^2),\\
\gamma_1\gamma_2^2\omega^2\int_0^{T_0}\xi_1(t)\xi_2(t)\,dt\le\frac{\gamma_1\gamma_2^3\omega^2}{2}\vert\xi_1\vert_{L^2(0,T_0)}^2+\frac{\gamma_1\gamma_2\omega^2}{2}\vert\xi_2\vert_{L^2(0,T_0)}^2,\\
\gamma_1\gamma_2\omega^2\xi_1^0\xi_2^0\le \displaystyle \frac{\gamma_1\gamma_2\omega^2}{2}(\vert \xi_1^0\vert^2+\vert \xi_2^0\vert^2),\\
\gamma_1\omega^2\xi_1(T_0)v^1_x(T_0,0)\le2\gamma_1\omega^4\vert\xi_1(T_0)\vert^2+ \frac{\gamma_1}{4}\vert v^1_x(T_0,0)\vert^2
\end{gather*}
together with \eqref{eq26*} and \eqref{eq45}, we finally get the existence of a positive constant $C_{\omega,\gamma_1,\gamma_2}$ which only depends on $\omega$, $\gamma_1$ and $\gamma_2$ (we underline the fact that $C_{\omega,\gamma_1,\gamma_2}$ does not depend on the time $T_0$) such that  \eqref{eq45} becomes
\begin{eqnarray*}
\begin{array}{l}
\label{eq46}
\displaystyle \frac{1}{2}(\vert v^2_{t}(T_0,.)\vert^2_{L^2(0,1)}+\vert v^2_{x}(T_0,.)\vert^2_{L^2(0,1)}+\gamma_1\omega^4\vert\xi_1(T_0)\vert^2)+\frac{\gamma_1}{4}\vert v^1_x(T_0,0)\vert^2\\
\displaystyle\phantom{gglllllkkkkg}
+\frac{\gamma_1\gamma_2\omega^2}{2}\int_0^{T_0}\vert \xi_2(t)\vert^2\,dt+\gamma_1\omega^2\vert \xi_2(T_0)\vert^2\le \\
\displaystyle \phantom{ggkkkkkkkkklllllllllllllkkg} C_{\omega,\gamma_1,\gamma_2}(\vert q^0\vert^2_{H^2(0,1)}+\vert q^1\vert^2_{H^1(0,1)}+\vert \xi_1^0\vert^2+\vert \xi_2^0\vert^2).
\end{array}
\end{eqnarray*}
This ends the proof of the energy estimates. An immediate consequence of these energy estimates is the uniqueness of the solution  $(v^1,v^2,\xi_1,\xi_2,\xi_3)\in C([0,T_0];L^2(0,1))\times C([0,T_0];H^{-1}(0,1))\times H^1(0,T_0)\times L^2(0,T_0)\times H^2(0,T_0)$ to \eqref{eq65}-\eqref{eq65bis}. Another consequence is that this solution is actually more regular, more precisely, $$(v^1,v^2)\in C([0,T_0];H^2(0,1)\cap H^1_r(0,1))\times C([0,T_0];H^{1}_r(0,1)).$$ Finally, the fact that the right members of both estimates \eqref{eq26*} and \eqref{eq26bis*} do not depend on the time $T_0$ ensures that the solution to \eqref{eq65} does not blow up at time $t=T_0$ and can also be extended at the right of this time. In other words, the solution to \eqref{eq65} exists, for any time $T>0$. In the same way, we obtain the same result for the system backward (system \eqref{eq65*}-\eqref{eq65bis*}).\\
Existence and uniqueness of a solution $(v^1,v^2,\xi_1,\xi_2,\xi_3)\in L_{loc}^{\infty}((0,+\infty); H^2(0,1)\cap H^1_r(0,1))\times L_{loc}^{\infty}((0,+\infty);H^{1}_r(0,1))\times H^1_{loc}(0,+\infty))\times L^2_{loc}(0,+\infty))\times H^2_{loc}(0,+\infty))$ to system \eqref{eqerror}-\eqref{eqerrorbis} 
immediately follows.\\
Let us now prove that energy inequalities are global in time, i.e., let us now prove \eqref{eq26} and \eqref{eq26bis}.

%%%%%%%%%%%%%%%%%%%%%%%%%%%%%%%%%%%%%%%%%%%%%%%%%%%%%%%%%%%%%%%%%%%%%%%%%%%%%%%%%%%%%%%%%%%%%%%%

\subsection{Proof of energy inequalities}
\underline{Step 1: proof of \eqref{eq26}.}\vspace{0,1cm}

Let $t\in (0,+\infty)$. We assume that there exists $k_0\ge0$ such that $2k_0T\le t\le (2k_0+1)T$. (The proof of \eqref{eq26} in the case were $(2k_0+1)T\le t\le (2k_0+2)T$ readily follows and we omit it.) 
One can compute easily, using \eqref{eqerror} and \eqref{eqerrorbis},
\begin{eqnarray*}
\begin{array}{rcl}
\label{equation1}
\displaystyle\int_0^t\int_0^1 v_t^2v^2\,dx\,dt&=&\displaystyle\sum\limits_{k=0}^{k_0-1}(\int_{2kT}^{(2k+1)T}\int_0^1v^1_{xx}v^2\,dx\,dt-\int_{(2k+1)T}^{(2k+2)T}\int_0^1v^1_{xx}v^2\,dx\,dt)\\&&\displaystyle+\int_{2k_0T}^t\int_0^1 v^1_{xx}v^2\,dx\,dt\vspace{0,1cm}\\
&=&\displaystyle \sum\limits_{k=0}^{k_0-1}(\int_{2kT}^{(2k+1)T}v_{xx}^1v^1_t\,dx\,dt+\int_{(2k+1)T}^{(2k+2)T}\int_0^1v_{xx}^1v^1_t\,dx\,dt\\&&\displaystyle+\int_{2k_0T}^t\int_0^1v_{xx}^1v^1_t\,dx\,dt\vspace{0,1cm}\\
&=&\displaystyle -\frac{1}{2}\int_0^t\int_0^1\partial_t(\vert v^1_x\vert^2)\,dx\,dt- \sum\limits_{k=0}^{k_0-1}\int_{2kT}^{(2k+1)T}(\dot{\xi}_2+\omega^2\xi_1)\gamma_1\xi_2\,dt\\
&&\displaystyle - \sum\limits_{k=0}^{k_0-1}\int_{(2k+1)T}^{(2k+2)T}(-\dot{\xi}_2+\omega^2\xi_1)(-\gamma_1\xi_2)\,dt)
 -\int_{2k_0T}^t(\dot{\xi}_2+\omega^2\xi_1)\gamma_1\xi_2\,dt\vspace{0,1cm}\\
 &=& \displaystyle  -\frac{1}{2}\int_0^t\int_0^1\partial_t(\vert v^1_x\vert^2)\,dx\,dt \\
 &&\phantom{ttkkkkkkkttttt}\displaystyle+\sum\limits_{k=0}^{k_0-1}\int_{2kT}^{(2k+1)T}(-\gamma_1\xi_2\dot{\xi}_2-\gamma_1\omega^2\xi_1(\dot{\xi}_1+\gamma_2\xi_1))\,dt\vspace{0,1cm}\\
 &&  \phantom{ttkkkkkkkttttt} \displaystyle + \sum\limits_{k=0}^{k_0-1}\int_{(2k+1)T}^{(2k+2)T}(-\gamma_1\xi_2\dot{\xi}_2+\gamma_1\omega^2\xi_1(-\dot{\xi}_1-\gamma_2\xi_1))\,dt\vspace{0,1cm}\\
 &&\phantom{ttkkkkkkkkkkkkkkkkttttt}  \displaystyle +\int_{2k_0T}^t(-\gamma_1\xi_2\dot{\xi}_2-\gamma_1\omega^2\xi_1(\dot{\xi}_1+\gamma_2\xi_1))\,dt.
\end{array}
\end{eqnarray*}
The end of the proof of \eqref{eq26} follows easily. We now focus on \eqref{eq26bis}.
\newpage
\underline{Step 2: proof of \eqref{eq26bis}.}\vspace{0,1cm}

We still assume that we are in the case where $2k_0T\le t\le (2k_0+1)T$. Then,
\begin{eqnarray*}
\begin{array}{rcl}
\label{equation2}
\displaystyle \int_0^t\int_0^1 v_{tt}^2v^2_t\,dx\,ds&=&\displaystyle\sum\limits_{k=0}^{k_0-1}(\int_{2kT}^{(2k+1)T}\int_0^1v^1_{xxt}v^2_t\,dx\,ds-\int_{(2k+1)T}^{(2k+2)T}\int_0^1v^1_{xxt}v^2_t\,dx\,ds)\\&&\displaystyle+\int_{2k_0T}^t\int_0^1 v^1_{xxt}v^2_t\,dx\,ds\vspace{0,1cm}\\
&=&\displaystyle -\frac{1}{2}\int_0^t\int_0^1\partial_t(\vert v^1_{xt}\vert^2)\,dx\,ds -\sum\limits_{k=0}^{k_0-1}\int_{2kT}^{(2k+1)T}v^1_{xt}(s,0)v^1_{tt}(s,0)\,ds\\
&&\displaystyle-\sum\limits_{k=0}^{k_0-1}\int_{(2k+1)T}^{(2k+2)T}v^1_{xt}(s,0)v^1_{tt}(s,0)\,ds-\int_{2k_0T}^tv^1_{xt}(s,0)v^1_{tt}(s,0)\,ds\vspace{0,1cm}\\
&=&\displaystyle -\frac{1}{2}\int_0^t\int_0^1\partial_t(\vert v^2_{x}\vert^2)\,dx\,ds\vspace{0,1cm}\\
&&\phantom{tttt}\displaystyle -\sum\limits_{k=0}^{k_0-1}\int_{2kT}^{(2k+1)T}v^1_{xt}(s,0)\gamma_1(-\omega^2\xi_1+v^1_x(s,0))\,ds\vspace{0,1cm}\\
&&\phantom{tttt}\displaystyle-\sum\limits_{k=0}^{k_0-1}\int_{(2k+1)T}^{(2k+2)T}v^1_{xt}(s,0)\gamma_1(-\omega^2\xi_1+v^1_x(s,0))\,ds\vspace{0,1cm}\\
&&\phantom{ttttjjjjjj}\displaystyle-\int_{2k_0T}^tv^1_{xt}(s,0)\gamma_1(-\omega^2\xi_1+v^1_x(s,0))\,ds.
\end{array}
\end{eqnarray*}
This rewrites
\begin{eqnarray*}
\begin{array}{rcl}
\label{equation3}
\displaystyle\frac{1}{2}\int_0^t\int_0^1(\partial_t(\vert v^2_t\vert^2)+\partial_t(\vert v^2_{x}\vert^2))\,dx\,ds+\frac{\gamma_1}{2}\int_0^t\partial_t(\vert v^1_x(s,0)\vert^2)\,ds&=&\displaystyle\gamma_1\omega^2\int_0^tv^1_{xt}(s,0)\xi_1(s)\,ds.
\end{array}
\end{eqnarray*}
Now, we still proceed as in the proof of \eqref{eq26bis*}
\begin{eqnarray}
\begin{array}{rcl}
\label{equation4}
\displaystyle \gamma_1\omega^2\int_0^tv^1_{xt}(s,0)\xi_1(s)\,ds&=&\displaystyle-\gamma_1\omega^2\int_0^{t}\dot{\xi_1}(s)v^1_x(s,0)\,ds+\gamma_1\omega^2[\xi_1(s)v^1_x(s,0)]^{t}_0\vspace{0,1cm}\\
&=&\displaystyle-\gamma_1\omega^2\sum\limits_{k=0}^{k_0-1}\int_{2kT}^{(2k+1)T}\dot{\xi_1}(s)(\dot{\xi_2}(s)+\omega^2\xi_1(s))\,ds\vspace{0,1cm}\\
&&\displaystyle-\gamma_1\omega^2\sum\limits_{k=0}^{k_0-1}\int_{(2k+1)T}^{(2k+2)T}\dot{\xi_1}(s)(-\dot{\xi_2}(s)+\omega^2\xi_1(s))\,ds\vspace{0,1cm}\\
&&\displaystyle-\gamma_1\omega^2\int_{2k_0T}^t\dot{\xi_1}(s)(\dot{\xi_2}(s)+\omega^2\xi_1(s))\,ds+\gamma_1\omega^2[\xi_1(s)v^1_x(s,0)]^{t}_0\vspace{0,1cm}\\
&=&\displaystyle-\frac{\gamma_1\omega^4}{2}(\vert\xi_1({t})\vert^2-\vert\xi_1^0\vert^2)+\gamma_1\omega^2\sum\limits_{k=0}^{k_0-1}\int_{2kT}^{(2k+1)T}\ddot{\xi}_1(s)\xi_2(s)\,ds\vspace{0,1cm}\\
&&\displaystyle-\gamma_1\omega^2\sum\limits_{k=0}^{k_0-1}\int_{(2k+1)T}^{(2k+2)T}\ddot{\xi}_1(s)\xi_2(s)\,ds+\gamma_1\omega^2\int_{2k_0T}^t\ddot{\xi}_1(s)\xi_2(s)\,ds\vspace{0,1cm}\\
&&\displaystyle-\gamma_1\omega^2\sum\limits_{k=0}^{k_0-1}[\dot{\xi}_1(s)\xi_2(s)]_{2kT}^{(2k+1)T}+\gamma_1\omega^2\sum\limits_{k=0}^{k_0-1}[\dot{\xi}_1(s)\xi_2(s)]_{(2k+1)T}^{(2k+2)T}\vspace{0,1cm}\\
&&\displaystyle-\gamma_1\omega^2[\dot{\xi}_1(s)\xi_2(s)]_{2k_0T}^{t}+\gamma_1\omega^2[\xi_1(s)v^1_x(s,0)]^{t}_0.
\end{array}
\end{eqnarray}
Let $A:=\displaystyle\gamma_1\omega^2(-\sum\limits_{k=0}^{k_0-1}[\dot{\xi}_1(s)\xi_2(s)]_{2kT}^{(2k+1)T}+\sum\limits_{k=0}^{k_0-1}[\dot{\xi}_1(s)\xi_2(s)]_{(2k+1)T}^{(2k+2)T}\displaystyle-[\dot{\xi}_1(s)\xi_2(s)]_{2k_0T}^{t}).$

\noindent One can compute
\begin{eqnarray*}
\begin{array}{rcl}
\label{equation5}
A&=&\gamma_1\omega^2\Big(-\sum\limits_{k=0}^{k_0-1}[\xi^2_2(s)-\gamma_2\xi_1(s)\xi_2(s)]_{2kT}^{(2k+1)T}+\sum\limits_{k=0}^{k_0-1}[-\xi_2^2(s)-\gamma_2\xi_1(s)\xi_2(s)]_{(2k+1)T}^{(2k+2)T}\\
&&\phantom{hhhhhhl:::::::::::::::::::::::::::lllllllllllllllllllllhh}-[\xi^2_2(s)-\gamma_2\xi_1(s)\xi_2(s)]_{2k_0T}^{t}\Big)\\
&=&\gamma_1\omega^2\Big(-\xi_2^2(t)+\xi_2^2(0)+\gamma_2\sum\limits_{k=0}^{k_0-1}(\xi_1((2k+1)T)\xi_2((2k+1)T)-\xi_1(2kT)\xi_2(2kT))\\
&&\displaystyle -\gamma_2\sum\limits_{k=0}^{k_0-1}(\xi_1((2k+2)T)\xi_2((2k+2)T)-\xi_1((2k+1)T)\xi_2((2k+1)T))\vspace{0,1cm}\\
&&+\gamma_2(\xi_1(t)\xi_2(t)-\xi_1(2k_0T)\xi_2(2k_0T))\Big)\vspace{0,2cm}\\
&=& \gamma_1\omega^2\Big(-\xi_2^2(t)+\xi_2^2(0)\vspace{0,1cm}\\
&&+\gamma_2\sum\limits_{k=0}^{k_0-1}(\xi_1((2k+1)T)(\dot{\xi}_1((2k+1)T)+\gamma_2\xi_1((2k+1)T))-\xi_1(2kT)(\dot{\xi}_1(2kT)+\gamma_2\xi_1(2kT)))\\&&
\displaystyle -\gamma_2\sum\limits_{k=0}^{k_0-1}(\xi_1((2k+2)T)(-\dot{\xi}_1((2k+2)T)-\gamma_2\xi_1((2k+2)T))\vspace{0,1cm}\\
&&\phantom{ttttttttttttttttttt}-\xi_1((2k+1)T)(-\dot{\xi}_1((2k+1)T)-\gamma_2\xi_1((2k+1)T)))\vspace{0,1cm}\\
&&+\gamma_2(\xi_1(t)(\dot{\xi}_1(t)+\gamma_2\xi_1(t))-\xi_1(2k_0T)(\dot{\xi}_1(2k_0T)+\gamma_2\xi_1(2k_0T)))\Big)\vspace{0,2cm}\\
&=&\gamma_1\omega^2(-\xi_2^2(t)+\xi_2^0+\gamma_2\xi_1(t)\xi_2(t)-\gamma_2\xi_1^0\xi_2^0).
\end{array}
\end{eqnarray*}
Thus, \eqref{equation4} becomes
\begin{eqnarray*}
\begin{array}{rcl}
\label{equation4finale}
\displaystyle \gamma_1\omega^2\int_0^tv^1_{xt}(s,0)\xi_1(s)\,ds&=&\displaystyle-\frac{\gamma_1\omega^4}{2}(\vert\xi_1({t})\vert^2-\vert\xi_1^0\vert^2)+\gamma_1\omega^2\sum\limits_{k=0}^{k_0-1}\int_{2kT}^{(2k+1)T}\ddot{\xi}_1(s)\xi_2(s)\,ds\vspace{0,1cm}\\
&&\displaystyle-\gamma_1\omega^2\sum\limits_{k=0}^{k_0-1}\int_{(2k+1)T}^{(2k+2)T}\ddot{\xi}_1(s)\xi_2(s)\,ds+\gamma_1\omega^2\int_{2k_0T}^t\ddot{\xi}_1(s)\xi_2(s)\,ds\vspace{0,1cm}\\
&&+\gamma_1\omega^2(-\xi_2^2(t)+\xi_2^0+\gamma_2\xi_1(t)\xi_2(t)-\gamma_2\xi_1^0\xi_2^0)+\gamma_1\omega^2[\xi_1(s)v^1_x(s,0)]^{t}_0.
\end{array}
\end{eqnarray*}
Now, for any $0\le k\le k_0$
\begin{gather*}
\label{equation6}
\displaystyle \gamma_1\omega^2\int_{2kT}^{(2k+1)T}\ddot{\xi}_1(s)\xi_2(s)\,ds=\displaystyle\gamma_1\omega^2\int_{2kT}^{(2k+1)T}(\dot{\xi}_2(s)-\gamma_2\xi_2(s)+\gamma_2^2\xi_1(s))\xi_2(s)\,ds,\vspace{0,1cm}\\
\label{equation7}
\displaystyle -\gamma_1\omega^2\int_{(2k+1)T}^{(2k+2)T}\ddot{\xi}_1(s)\xi_2(s)\,ds=\displaystyle-\gamma_1\omega^2\int_{2kT}^{(2k+1)T}(-\dot{\xi}_2(s)+\gamma_2\xi_2(s)+\gamma_2^2\xi_1(s))\xi_2(s)\,ds,\vspace{0,1cm}\\
\label{equation8}
\gamma_1\omega^2\int_{2k_0T}^t\ddot{\xi}_1(s)\xi_2(s)\,ds=\displaystyle \gamma_1\omega^2\int_{2k_0T}^t(\dot{\xi}_2(s)-\gamma_2\xi_2(s)+\gamma_2^2\xi_1(s))\xi_2(s)\,ds,
\end{gather*}
thus
\begin{eqnarray*}
\begin{array}{l}
\label{equation9}
\displaystyle\gamma_1\omega^2(\sum\limits_{k=0}^{k_0-1}\int_{2kT}^{(2k+1)T}\ddot{\xi}_1(s)\xi_2(s)\,ds\displaystyle-\sum\limits_{k=0}^{k_0-1}\int_{(2k+1)T}^{(2k+2)T}\ddot{\xi}_1(s)\xi_2(s)\,ds+\int_{2k_0T}^t\ddot{\xi}_1(s)\xi_2(s)\,ds)=\vspace{0,1cm}\\
\displaystyle\gamma_1\omega^2\Big(\frac{1}{2}\int_0^t\partial_t(\vert \xi_2(s)\vert^2\,ds-\gamma_2\int_0^t\vert \xi_2(s)\vert^2\,ds+\gamma_2^2\sum\limits_{k=0}^{k_0-1}\int_{2kT}^{(2k+1)T}\xi_1(s)(\dot{\xi}_1(s)+\gamma_2\xi_1(s))\,ds\vspace{0,1cm}\\
\displaystyle-\gamma_2^2\sum\limits_{k=0}^{k_0-1}\int_{(2k+1)T}^{(2k+2)T}\xi_1(s)(-\dot{\xi}_1(s)-\gamma_2\xi_1(s))\,ds+\gamma_2^2\int_{2k_0T}^t\xi_1(s)(\dot{\xi}_1(s)+\gamma_2\xi_1(s))\,ds\Big).
\end{array}
\end{eqnarray*}
So, finally,
\begin{eqnarray*}
\begin{array}{l}
\label{equation10}
\displaystyle\gamma_1\omega^2\Big(\sum\limits_{k=0}^{k_0-1}\int_{2kT}^{(2k+1)T}\ddot{\xi}_1(s)\xi_2(s)\,ds\displaystyle-\sum\limits_{k=0}^{k_0-1}\int_{(2k+1)T}^{(2k+2)T}\ddot{\xi}_1(s)\xi_2(s)\,ds+\int_{2k_0T}^t\ddot{\xi}_1(s)\xi_2(s)\,ds\Big)=\vspace{0,1cm}\\
\displaystyle\gamma_1\omega^2\Big(\frac{1}{2}\int_0^t\partial_t(\vert \xi_2(s)\vert^2\,ds-\gamma_2\int_0^t\vert \xi_2(s)\vert^2\,ds+\frac{1}{2}\gamma_2^2\int_0^t\partial_t(\vert \xi_1(s)\vert^2)\,ds+\gamma_2^3\int_0^t\vert \xi_1(s)\vert^2\,ds\Big).
\end{array}
\end{eqnarray*}
Following now the proof of \eqref{eq26bis*}, one easily obtains \eqref{eq26}.
%-----------------------------------------------------------------------------------------
%---------------------------------
%-----------------------------------------------------------------------------------------
%---------------------------------%              ----------------------                       ------------------------------------
The proof of Proposition \ref{prop2} is now complete. \cqfd
%----------------------------------------------------------------------------------------------------------------
%----------------------------------------------------------------------------------------------------------------

\def\cprime{$'$} \newcommand{\SortNoop}[1]{}

\end{document}